\documentclass[lefttitle, final,5p,times,twocolumn]{elsarticle}

\usepackage{stfloats}
\usepackage[utf8]{inputenc}
\usepackage{etoolbox}
\usepackage[T1]{fontenc}
\usepackage{lmodern}
\usepackage[english]{babel}
\usepackage{amsmath,amssymb, amsthm,bm}
\usepackage{newtxtext,newtxmath}
\usepackage{graphicx}
\usepackage{algorithm}
\usepackage{tabularx}
\usepackage{algpseudocode}
\usepackage{float}
\usepackage{lineno}
\usepackage{booktabs}
\usepackage{caption}
\usepackage{subcaption}
\usepackage{enumitem}
\usepackage{xcolor}
\usepackage{balance}
\usepackage{url}
\usepackage{placeins}
\usepackage{verbatim}

\usepackage{hyperref}
\usepackage[capitalize]{cleveref}

\abstracttitle{}

\DeclareMathOperator{\N}{\mathbb{N}}

\journal{Applied Mathematics and Computation}

\begin{document}

\theoremstyle{plain}
\newtheorem{thm}{Theorem}[section]
\newtheorem{proposition}{Proposition}[section]
\newtheorem{lem}[thm]{Lemma}
\newtheorem{cor}[thm]{Corollary}
\newtheorem{theorem}{Theorem}[section]
\theoremstyle{definition}
\newtheorem{defn}[thm]{Definition}
\newtheorem*{notation}{Notation}
\newtheorem{example}[thm]{Example}
\newtheorem{quest}[thm]{Question}
\newtheorem{conj}[thm]{Conjecture}

\theoremstyle{remark}
\newtheorem{rem}[thm]{Remark}

\renewcommand{\algorithmicrequire}{\textbf{Input:}}
\renewcommand{\algorithmicensure}{\textbf{Output:}}

\makeatletter
\def\ps@pprintTitle{%
 \let\@oddhead\@empty
 \let\@evenhead\@empty
 \def\@oddfoot{\centerline{\thepage}}%
 \let\@evenfoot\@oddfoot}
\makeatother
\patchcmd{\abstract}{\null\vfil}{}{}{}
\patchcmd{\endabstract}{\par}{}{}{}
\journal{Physica D: Nonlinear Phenomena}


\begin{frontmatter}


\title{Numerical Transitivity and Numerical Leo Properties for Lorenz Maps with Applications to Courbage-Nekorkin-Vdovin Neuron Model}

\author[inst1]{Rudrakshala Kavya Sri}
\author[inst2]{Piotr Bart{\l}omiejczyk}
\author[inst1]{Sishu Shankar Muni\corref{cor1}}
\cortext[cor1]{Corresponding author}
\ead{sishushankarmuni@gmail.com}

\address[inst1]{School of Digital Sciences, Digital University Kerala, 695317, Pallipuram, India}
\address[inst2]{Faculty of Applied Physics and Mathematics, Gabriela Narutowicza 11/12, 80-233 Gdańsk, Poland}

\begin{abstract}
\small
\noindent
\parbox[t]{0.28\textwidth}{
  {\large{ARTICLE INFO}}\\\vspace{2pt}%
   \noindent\hrule
   \vspace{4pt}  
   
  secondary 92B20; 37M10\\%
  
  Locally eventually onto\\
  Chaos\\
  Lorenz map\\
  Neuron model
}
\hspace{0.04\textwidth}
\parbox[t]{0.67\textwidth}{
  {\large{ABSTRACT}}\\\vspace{2pt}%
  
  \noindent\hrule
  \vspace{4pt}
  This research investigates the dynamic behavior of one-dimensional discrete systems using two computational algorithms: the \emph{numerical transitivity} and the \emph{numerical locally eventually onto} (LEO) tests. Both algorithms are systematically applied to a variety of interval maps, including classical examples such as $\beta$-transformations and expanding Lorenz maps, in order to assess and characterize their chaotic dynamics. We perform a detailed comparison of the two methods in terms of accuracy, computational efficiency, and their sensitivity in detecting transitions between regular and chaotic regimes.
  Particular emphasis is placed on the Courbage–Nekorkin–Vdovin (CNV) model of a single neuron, known for its rich, spiking-like dynamics and its mathematical reducibility to Lorenz-type maps. By analyzing both the piecewise linear and nonlinear versions of the CNV model, we illustrate how the proposed numerical tests reliably capture qualitative changes in the system’s dynamics,  focusing on the onset of chaos and chaotic regimes. The results highlight the practical potential of these numerical approaches as diagnostic tools for studying complex dynamical systems arising in nonlinear science and mathematical neuroscience.}
\end{abstract}

\end{frontmatter}
\section{Introduction}

This paper presents two numerical 
algorithms—\emph{numerical transitivity} 
and \emph{numerical LEO}
(LEO abbreviates \emph{locally eventually onto})—for analyzing 
dynamical behavior of one-dimensional maps, 
particularly Lorenz maps and their types. 
These algorithms provide computational methods for detecting 
key features related to chaotic dynamics: 
topological transitivity and topological LEO. 
We compare operation and results of both algorithms 
and apply them to various classes of Lorenz maps. 
As a central application, we examine their effectiveness 
on the Courbage–Nekorkin–Vdovin (CNV for short) neuron model, 
both in its piecewise linear and nonlinear versions, 
to study how these maps 
behave under different parameter settings.

In the early 1960s, Edward Lorenz demonstrated 
that small variations in initial conditions can result 
in drastically different outcomes, and thus giving rise to 
\emph{Lorenz maps} \cite{lorenz1963deterministic}. 
Lorenz maps illustrate the butterfly effect
in the sense that they demonstrate how small variations in initial conditions 
can result in substantial differences over time, which is crucial 
in grasping chaotic systems and elucidating their complex, 
unpredictable movements within recognizable 
patterns \cite{analysis_of_chaotic_systems, palmer2014real, shen2022three}. 
A little bit earlier, in the late 1950s,
$\beta$-\emph{transformations} 
were developed in the number theory 
to more accurately represent real numbers
see~\cite{parry1960,renyi1957}. 
Both forms have expansion and are characterized by discontinuities, 
making them perfect for modeling sensitive and unpredictable systems. 
In our study, these maps offer a simple framework for the analysis of 
the dynamics of the CNV neuron model that possesses similar features 
of complexity and sensitivity.

Recall that, in its simplest form,
$\beta$-transformations are mathematical functions 
that stretch a number by a specific factor \(\beta\),
retaining only its decimal part. 
These transformations are important for number representation, 
chaotic dynamics, and cryptography, 
aiding in improved precision and data compression. 
Chaotic behavior of $\beta$-transformations
helps reveal repeating 
and self-similar patterns in fractals and uncover how states 
evolve over time in dynamical systems, making them useful 
for studying complex, unpredictable patterns 
\cite{beta_fractal,beta_type}. 
These transformations play 
a crucial role in cryptography, ergodic optimization, 
ensuring secure communications, and are also used in the modeling 
of specific types of random processes
\cite{beta_applications, 
beta_in_distribution,beta_transformation}. 

\emph{Lorenz maps} were introduced in the study of the Lorenz system 
and attractor and to this they owe their name
\cite{lorenz1963deterministic,guckenheimer1979stability,
classic_lorenz_attractor}.
The classical Lorenz attractor demonstrates \emph{mixing}
behavior \cite{lorenz_attractor_as_mixing}. 
Mixing is the process where, over time, an initial condition 
in a system blends with others, making the future behavior 
of the system independent of its starting point. 
This lack of memory of the past is a key characteristic 
of chaos\cite{fractional_lorenzmaps}.
It highlights the relevant features of the Lorenz attractor 
that can be studied using one-dimensional dynamical systems, 
particularly interval maps with discontinuities. 
Furthermore, specific classes of geometric Lorenz 
flows, introduced by Afraimovich, Bykov, Shilnikov
\cite{afraimovich1977origin} and Guckenheimer, 
Williams \cite{guckenheimer1979stability}, 
can also exhibit mixing under certain conditions. 
Let us emphasize that transitivity and LEO are 
also properties of the generalized mixing type.

\emph{Expanding Lorenz maps} provide natural and visual examples 
of chaotic behavior through mixing, 
where it can take points in the domain and distribute them across 
the space so that they cover more space \cite{mixing}. 
Topological mixing inhibits patterns, 
so that the points continue to develop towards more randomness. 
Strong mixing in the system tends to forget its initial condition over 
time\cite{lorenz_attractor_as_mixing}.  It makes expanding Lorenz maps 
useful in applications where simple versions of one-dimensional 
nonlinearity are emerging.

\emph{Chaos} in interval maps can be defined in several ways. 
\emph{Li-Yorke chaos} occurs when some point pairs act in a random manner, 
occasionally approaching one another and then moving apart. 
If there are a large number of such points, then the map is said 
to be Li-Yorke chaotic, a fairly weak form of chaos that occurs in most systems. 
\emph{Entropy chaos} measures the complexity and the high level of unpredictability 
of the dynamics of the system. A positive entropy indicates strong chaos, 
with many possible behaviors. \emph{Devaney chaos} is a stricter class and can be  
defined by the coexistence of sensitive dependence on initial conditions, 
topological transitivity, and a dense set of periodic 
points \cite{chaos_for_continuous_maps}.
In this paper we will mainly focus on Devaney chaos.
Although it is worth emphasizing that the functions 
we are considering have discontinuities.

\emph{Transitivity} in Lorenz maps is necessary for overall mixing to occur 
and is an aspect of chaotic motion. It implies that every point can come close 
to any other point in the long run. For one-dimensional Lorenz maps, 
random behaviors and transitivity can be induced by specific conditions. 
A transitive Lorenz map can be expected to be chaotic 
and have strongly sensitive dynamics
\cite{devaney1989definition, guckenheimer1979stability}. 
The determination of the parameter ranges for transitivity assists 
in differentiating between general mixing and other behaviors. 
These characteristics have real-world applications in fluid 
dynamics \cite{lorenz1963deterministic}, improving forecasts and 
sharpening mathematical models in climate modeling and turbulence 
studies \cite{transitivity_in_lorenz_maps, Dynamical_system}.

\emph{Neurons} are important in the nervous system to transmit information using 
electrical and chemical signals. In turn, neuron models such 
as the \emph{Courbage-Nekorkin-Vdovin model} 
\cite{courbage2007,courbage2010} 
are mathematical descriptions that 
allow neuroscience and applied mathematics 
to study neural dynamics and chaotic systems,
simplifying neuronal behavior in a discrete-time framework 
\cite{multilevel, cnv_model, cnv_model1}. 
Neuronal \emph{spiking} occurs when a neuron sends a rapid electrical signal, 
called an action potential, to communicate with other neurons. 
This action is effectively a short message for those neurons 
to interact \cite{spiking_neurons}. 
The CNV model can be used to reproduce features of neuronal activity kinetics, 
primarily with respect to various spike patterns (both regular and chaotic). 
For the CNV model, the numerical transitivity and LEO tests
help determine whether the neuron's membrane potential 
crosses its state space thoroughly (transitivity) or completely (LEO), 
providing deeper insights into chaotic behavior.

The CNV model was introduced by Courbage, Nekorkin, and Vdovin, 
as a minimal model of the dynamics of neurons firing, accomplished 
via one-dimensional piecewise continuous maps. 
Since the CNV model map is piecewise continuous and monotonic
(in a finite interval) 
with one critical discontinuity, it is formally equivalent 
to an expanding Lorenz map. In an effort to offer 
more flexibility to CNV analyses, 
the authors of~\cite{courbage2007,courbage2010}  
provided both piecewise linear (plCNV) 
and nonlinear (nlCNV) versions of the model. 
The plCNV model seems to be easier to analyze and is equivalent to
a $\beta$-transformation. On the other hand, the nlCNV model 
presents a more mathematically difficult challenge 
due to less obvious but equally complex dynamics. 
Both models provide a realistic and adjustable simulation environment 
for exploring shifts between regular spiking, bursting, 
and chaotic patterns that are seen in real neuronal behavior.

The organization of the paper is as follows. Section \ref{sec:prel} 
contains short preliminaries, 
including definitions of transitivity, periodic points, 
and Devaney chaos. Section \ref{sec:lorenz} discusses a number of classes 
of Lorenz maps and describes their dynamical properties. Sections \ref{sec:numtrans} 
and \ref{sec:numleo} provides a detailed exposition of 
the numerical transitivity and LEO algorithms,
which are the primary focus of our research,
and give their mathematical understandings. 
In turn, in Section \ref{sec:comparison}
we present detailed comparison
of numerical transitivity and LEO algorithms.
Section \ref{sec:cnv} discusses 
the CNV neuron model in both piecewise linear and nonlinear contexts. 
Section \ref{sec:transleocnv} presents the application of the transitivity 
and LEO algorithms to analyze the behavior of the CNV model 
across a range of parameters. Finally, Section \ref{sec:discussion} 
concludes the paper with a summary and disussion.
All figures in this paper have been prepared 
with MATLAB’s implementation of the algorithms presented here.

\section{Preliminaries}
\label{sec:prel}

Let $(X,d)$ be a metric space and $f\colon X \to X$ a map (not necessarily continuous).
Recall that $f^n$ denotes the $n$-fold composition of $f$ with itself and, if $x\in X$, 
then the orbit of $x$ under $f$ is the set 
\[
O(x):=\{f^n(x)\mid n\ge0\}.
\]

\subsection{Periodic points}

A point $x\in X$ such that $f(x)=x$ is \emph{periodic} (of period $n\in\N$, $n\ge1$), if
\[
f^n(x)=x \quad \text{and} \quad f^k(x)\neq x \quad \text{for } k=1,\dotsc,n-1.
\]
An orbit of a periodic point is called a \emph{periodic orbit}.

\subsection{Topological transitivity, topological mixing and sensitivity}
\label{subsec:trans}

Recall that $f$ is called
\begin{itemize}
    \item \emph{topologically transitive} if for every two nonempty open sets
    $U,V\subset X$ there exists $n\in\N$ such that 
    $f^n(U)\cap V\neq\emptyset$
    (otherwise we will call it \emph{topologically nontransitive}),
    \item \emph{topologically mixing} if for every two nonempty open sets 
    $U,V\subset X$ there exists $N\in\N$ such that for every 
    $n\ge N$ we have $f^n(U)\cap V\neq\emptyset$,
    \item \emph{sensitive} if there exists  $\delta>0$ such that
    for every $x\in X$ and every neighborhood 
    $U$ of $x$, there exists $y\in U$ and $n\in\N$ 
    such that $d(f^n(x),f^n(y))>\delta$.
\end{itemize}

Note that, by definition, topological mixing 
implies topological transitivity.

\subsection{Chaos}

One of the best-known definitions of chaos is 
the following due to R.\ L.\ Devaney (1989).

A function $f\colon X\to X$ is called 
\emph{chaotic in the sense of Devaney} on $X$ if 
\begin{enumerate}
    \item $f$ is topologically transitive,
    \item the set of periodic points of $f$ is dense in $X$,
    \item $f$ is sensitive.
\end{enumerate}

\section{Lorenz maps}
\label{sec:lorenz}

The basic class of maps studied in this paper are Lorenz maps. 
Authors consider different types of Lorenz maps,
for example, $\beta$-transformations (R\'enyi, Parry),
Lorenz-type maps (Afraimovich),
Lorenz-like maps (Alsed\`a, Misiurewicz)
or expanding Lorenz maps (Raith, Oprocha).
Lorenz maps appear in various applications
as the Lorenz system
and modeling a single neuron
(Courbage-Nekorkin-Vdovin model).
Here we will mainly concentrate on
three subclasses of Lorenz maps:
\emph{Lorenz-like maps}, \emph{expanding Lorenz maps}
and $\beta$-\emph{transformations},
which are defined below.

For simplicity of notation, we will formulate 
the definitions below for 
the unit interval $[0,1)$. However,
all definitions make sense and all results still hold 
if we replace the unit interval 
by $[a,b)$ for fixed $a$ and $b$
and use the linear change 
of variables, which is a topological conjugacy.

\subsection{Lorenz-like maps}

A \emph{Lorenz-like map} is a map $f$ of an interval
$[0,1)$ to itself, for which there exists a point
$c\in(0,1)$ such that 
\begin{itemize}
	\item $f$ is continuous and increasing
        (not necessarily strictly) 
        on $[0, c)$ and on $(c, 1)$,
	\item $\lim_{x\to c^-}f(x) = 1$ and 
	$\lim_{x\to c^+}f(x) = f(c) = 0$.
\end{itemize} 

In Figure~\ref{fig:lorenz-expanding-beta}, the left panel shows 
a Lorenz-like map where the function \(f\) meets 
the given conditions 
\[
f(x) =
\begin{cases}
2/5, & 0\leq x \leq 1/5, \\
2/5 + (12/5)\cdot(x-1/5), & 1/5<x<9/20,\\
4(x-9/20), & 9/20\leq x \leq 3/5,\\
3/5, & 3/5<x< 1. \\
\end{cases}
\]
Note that a Lorenz-like map can be constant 
on some subinterval 
of the domain. Such a situation cannot occur 
for the two classes of maps considered 
in the following subsections.
Finally, observe that if we consider a Lorenz-like map on a circle $S^1$ rather 
than on an interval $[0,1)$, the discontinuity at the point $c$ disappears, 
but the discontinuity at $x=0$ appears instead.

\subsection{Expanding Lorenz maps}

An \emph{expanding Lorenz map} is a map
$f\colon[0,1)\to[0,1)$ satisfying 
the following three conditions:
\begin{itemize}
    \item there is a \emph{critical point}
    $c\in(0, 1)$ such that $f$ is continuous 
    and strictly increasing on $[0, c)$ and 
    $(c, 1)$,
    \item $\lim_{x\to c^-}f(x)=f(c^-)=1$ and 
    $\lim_{x\to c^+}f(x)=f(c^+)=f(c)=0$,
    \item $f$ is differentiable for all points
    not belonging to a finite set 
    $F\subset[0, 1)$ and 
    \[
    \beta_f:=
    \inf{\{f'(x)\mid x\in[0,1)\setminus F\}}>1.
    \]
\end{itemize}

By definition, expanding Lorenz maps 
are Lorenz-like. The middle panel
of Figure~\ref{fig:lorenz-expanding-beta}
shows a nonlinear expanding Lorenz map,
and the right panel shows a constant slope (linear)
expanding Lorenz map.
The nonlinear map is given as
\[
f(x) = 
\begin{cases}
0.1 + 0.9 \cdot \dfrac{\exp\left(3x/2\right) - 1}{\exp\left(3c/2\right) - 1}, & \text{if } 0\le x < c, \\
0.9 \cdot \left(1 - \dfrac{\exp\left(3\left(1-x\right)/2\right) - 1}{\exp\left(3 \left(1-c\right)/2\right) - 1} \right), & \text{if } c\le x < 1,
\end{cases}
\]
where $c=0.45$.
We will see in Section~\ref{sec:cnv} that 
the 1D CNV neuron model is described by 
an expanding Lorenz map both 
in piecewise linear and nonlinear cases.

\subsection{$\beta$-transformations}

\begin{figure*}[!htb]
    \centering
    \includegraphics[scale=0.53, trim=0mm 2mm 0mm 5mm]{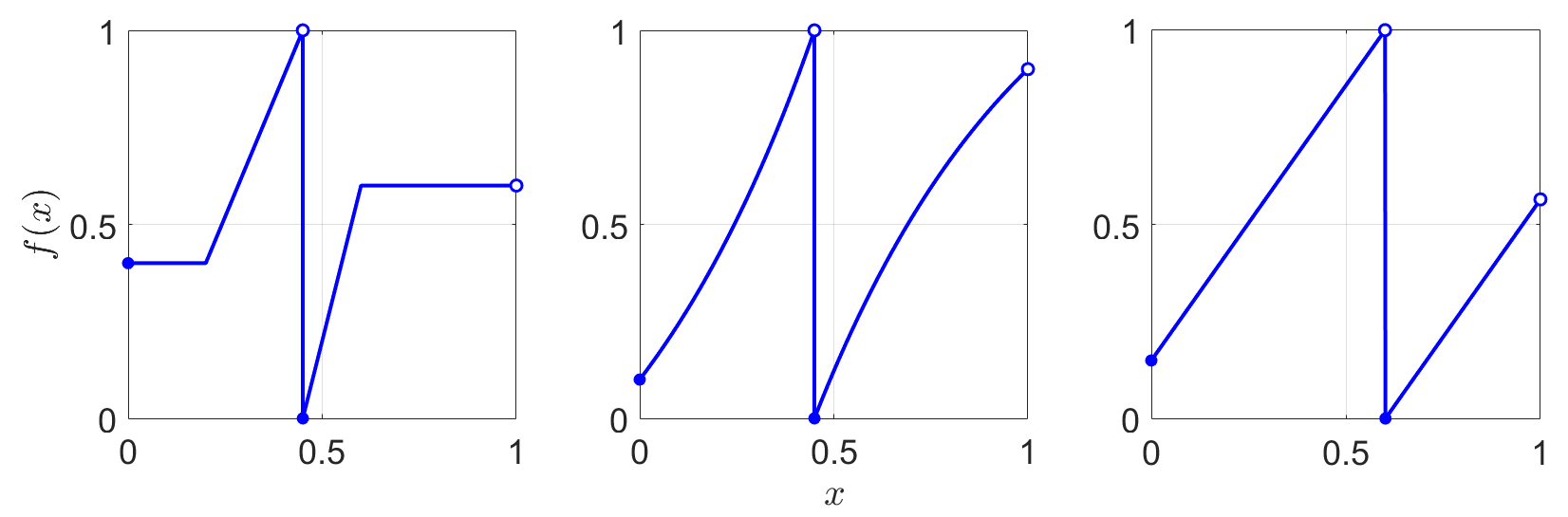}
    \caption{
    The left panel presents an example of a Lorenz-like map, the middle panel an expanding Lorenz map, and the right panel a $\beta$-transformation.
    }
    \label{fig:lorenz-expanding-beta}
\end{figure*}

Let $1<\beta\le2$, $\alpha\ge0$ and $\alpha+\beta\le2$.
The map $T\colon[0,1)\to[0,1)$ of the form
\[
T(x)=\beta x+\alpha\!\! \pmod 1
\]
is called a $\beta$-\emph{transformation}.
Note that every $\beta$-transformation is 
an expanding Lorenz map.
In particular, each $\beta$-transformation
has a unique point of discontinuity at
$x=(1-\alpha)/\beta$ and $T(x)=0$. 
The right panel of Figure~\ref{fig:lorenz-expanding-beta}
presents a $\beta$-transformation, i.e., 
a constant slope map given by the formula
$T(x)=\sqrt{2} x+3/20\! \pmod 1$.
The piecewise linear version of 
the 1D CNV model map is a $\beta$-transformation.
In Figure~\ref{fig:lorenz-expanding-beta}, we 
have considered the following formula 
to plot the $\beta$-transformation 
map that satisfies the above conditions:
\[
f(x) = \beta x+\alpha\!\! \pmod 1, \quad 
\text{where } \beta = \sqrt{2},\ \alpha = 0.15.
\]

\section{Numerical transitivity}
\label{sec:numtrans}

In this section we introduce the notion of \emph{numerical transitivity} and show its possible application in the study of Lorenz maps. It can be regarded as the numerical equivalent of the \emph{topological transitivity} property well known in the theory of dynamical systems (see Subsection~\ref{subsec:trans} for the definition). However, our algorithm of numerical transitivity is not directly based on the definition of topological transitivity.

Namely, the following result proved in~\cite[Prop.~1]{dense_orbit} provides a natural characterization of topological transitivity in terms of a dense orbit and is the main motivation for our algorithm of numerical transitivity.

\begin{proposition}\label{prop:dense}
A Lorenz-like map $f\colon [0,1)\to[0,1)$ is topologically transitive if and only if it has a dense orbit, i.e., there exists a point $x\in[0,1)$ such that every nonempty subinterval of $[0,1)$ contains an element of $O(x)$.
\end{proposition}

Why is the topological transitivity property so important? Because, for the class of expanding Lorenz maps, it is actually \emph{equivalent} to the presence of \emph{Devaney chaos} on the entire domain—that is, in the whole interval. This connection is formalized in the following theorem, proven in~\cite[Thm.~A.8]{cnv_model1}.

\begin{theorem}\label{thm:chaoslorenz}
Let $f:[0,1)\to [0,1)$ be an expanding Lorenz map. Then the following conditions are equivalent:
\begin{itemize}[leftmargin=2em]
    \item $f$ is topologically transitive,
    \item $f$ is chaotic in the sense of Devaney on $[0,1)$.
\end{itemize}
\end{theorem}

We emphasize that, assuming that our numerical transitivity test is valid (trustworthy), it is possible to determine from it whether an expanding Lorenz map is chaotic in the sense of Devaney on its whole domain or not. It is also worth recalling that a similar result holds for continuous maps on an interval (see~\cite{vellekoop}). However, it does not apply in our case because expanding Lorenz maps are not continuous by definition.

\subsection{Algorithm of numerical transitivity}

We are now ready to present, in points, the main idea of the numerical transitivity test.

\begin{enumerate}[label=\arabic*.]
    \item The algorithm tests whether a function \( f \) thoroughly explores its entire domain when iterated over time.
    \item It generates several sequences (trials) by repeatedly applying the function \( f \), starting from a random point within a given interval.
    \item For each sequence, it tracks how often the function’s outputs fall into evenly spaced subintervals (bins) of the domain.
    \item If, during any trial, all bins contain at least one value from the sequence, the function is considered to have the property of numerical transitivity.
    \item If no trial achieves complete coverage of all bins, the function is considered to have failed the test of numerical transitivity.
\end{enumerate}

\begin{algorithm}[H]
{\small
\caption{\textsc{Numerical Transitivity Test for Map $f$ on Interval $[a,b)$}}
\label{alg:TransTest}
\begin{algorithmic}[1]
\Function{NumTransTest}{$f, a, b$} \Comment{$f$ is a map defined on $[a,b)$ into itself}
    \State $N \gets 50000$ \Comment{number of iterations (length of time series)}
    \State $\text{\it{num\_trials}} \gets 5$ \Comment{number of trials}
    \State $k \gets 200$ \Comment{transient in time series}
    \State $\delta \gets (b-a)/1000$ \Comment{diameter of division}
    \State $\text{\it{ranges}} \gets a:\delta:b$ \Comment{division into bins}
    
    \For{$j \gets 1$ to $\text{\it{num\_trials}}$}
        \State $x \gets$ array of size $N$, initialized to $0$ \Comment{preallocation of time series}
        \State $x[1] \gets$ random value between $a$ and $b$ \Comment{random initial point}
        
        \For{$i \gets 2$ to $N$}
            \State $x[i] \gets f(x[i-1])$ \Comment{construction of time series}
        \EndFor
        \State $x \gets x[k:\text{end}]$ \Comment{removing transient from time series}
        \State $\text{\it{bincounts}} \gets$ histogram of $x$ using $\text{\it{ranges}}$ \Comment{counting elements in bins}
        \If{all values in $\text{\it{bincounts}} > 0$} \Comment{checking that all bins are nonempty}
            \State $\text{\it{logValue}} \gets$ \textbf{true} 
            \Comment{test passed}
            \State \Return $\text{\it{logValue}}$
            \Comment{end of program}
        \EndIf
    \EndFor
    \State $\text{\it{logValue}} \gets$ \textbf{false} \Comment{test failed}
    \State \Return $\text{\it{logValue}}$ \Comment{end of program}
\EndFunction
\end{algorithmic}}
\end{algorithm}

Based on the above guidelines, Algorithm~\ref{alg:TransTest} presents in pseudocode a test of the property of \emph{numerical transitivity} for an expanding Lorenz map $f\colon [a,b)\to[a,b)$. It is worth pointing out that this is just one of many possible implementations for checking the existence of a dense orbit. The algorithm has as
\begin{itemize}[leftmargin=2em]
    \item \textbf{input:} a map $f$ and two numbers $a$ and $b$ (endpoints of the domain of $f$), assuming that $f$ is an expanding Lorenz map,
    \item \textbf{output:} logical value $1$ (\textbf{true}) if $f$ is numerically transitive and $0$ (\textbf{false}) if $f$ is not.
\end{itemize}

\subsection{\texorpdfstring{Numerical transitivity simulations for $\beta$-transformations}{}}
\label{subsec:transsim}

In this subsection, we justify that the \emph{numerical transitivity} provides a computationally efficient and reliable substitute for the formal definition of \emph{topological transitivity}. Consequently, it can be used to determine with good approximation whether a discrete dynamical system given by an expanding Lorenz map is topologically transitive or not. 

It turns out that the results of the numerical transitivity test, on the one hand, confirm several well-known theoretical results for $\beta$-transformations, and on the other hand, make it possible to predict some new conjectures.

\begin{figure*}[!htb]
    \centering
    \includegraphics[scale=0.23, trim=5mm 5mm 0mm 0mm]{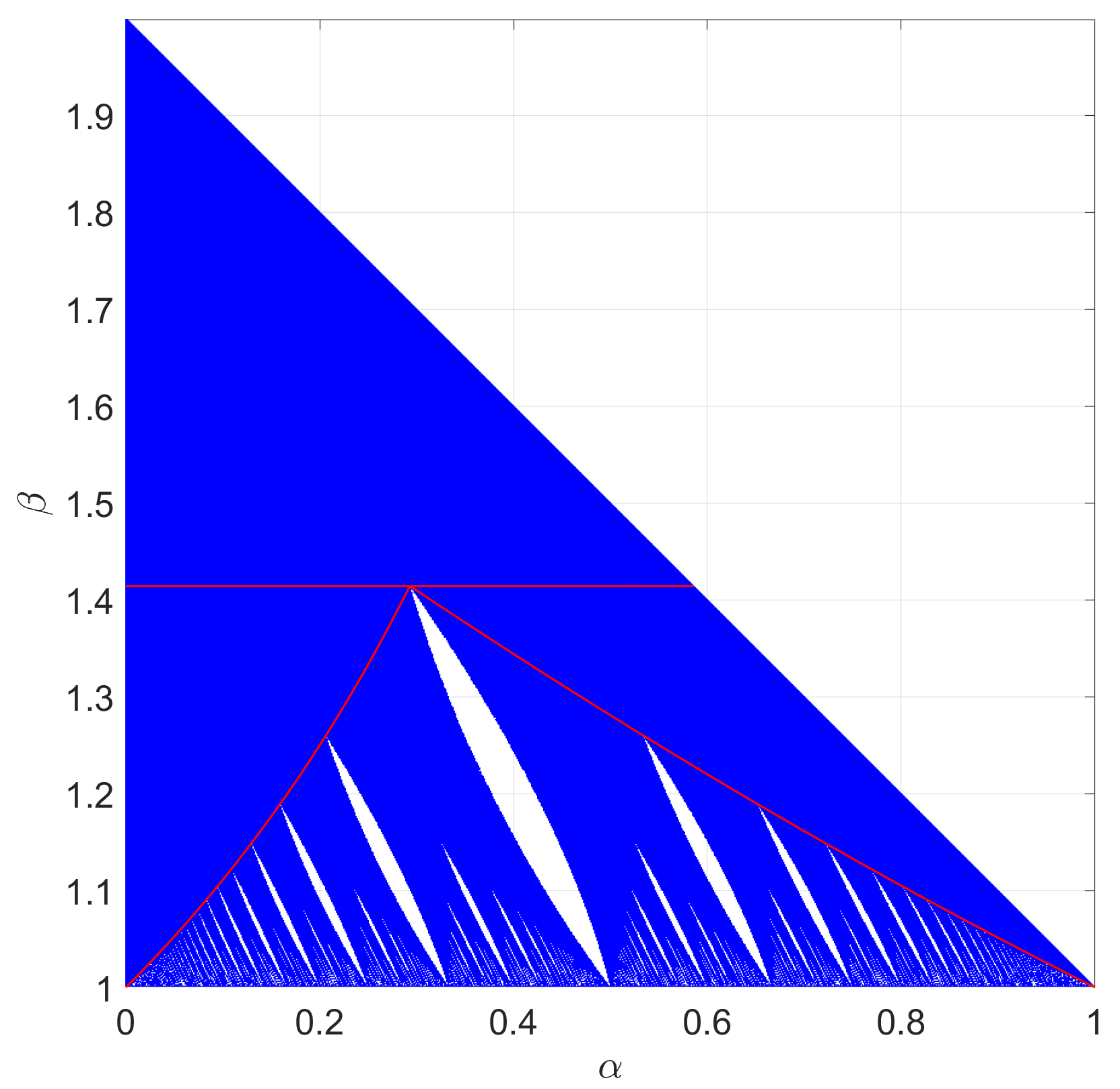}
    \caption{Results of the numerical transitivity test for classical $\beta$-transformations in the $\alpha$–$\beta$ parameter triangle $\mathcal{T}$. Horizontal red line: $\beta=\sqrt{2}$, left red line: $\alpha=1-1/\beta$, right red line: $\alpha=1+1/\beta-\beta$. Meshgrid = 800.}
    \label{fig:triangleTrans}
\end{figure*}

To begin with, consider the well-known right-angled triangle $\mathcal{T}$ of $\alpha$–$\beta$ parameters for $\beta$-transformations (see Fig.~\ref{fig:triangleTrans}, which is similar but not identical to illustrations in~\cite{mixing,palmerthesis}), given by the conditions:
\[
1<\beta\le2, \qquad \alpha\ge0, \qquad \alpha+\beta\le2.
\]
By definition, there is a natural bijective correspondence between all possible $\beta$-transformations and the points from the triangle $\mathcal{T}$. Moreover, the blue points in $\mathcal{T}$ correspond to the positive result of the numerical transitivity test for the respective $\beta$-transformation, while the white region (inside $\mathcal{T}$) corresponds to the negative result of the test.

With this in mind, and taking into account the fixed resolution of Fig.~\ref{fig:triangleTrans}, we can immediately, and graphically, confirm the following three statements proved in~\cite[Thm.~6.6, Thm.~4.6, Thm.~7.1]{mixing} concerning the topological transitivity property for $\beta$-transformations, i.e. maps $T\colon[0,1)\to[0,1)$ given by
\[
T(x)=\beta x+\alpha\!\! \pmod 1.
\]

\begin{enumerate}[label=\textbf{(\arabic*)}]
    \item Assume $n\ge2$ and $1<\beta<\sqrt[n]{2}$. If
    \[
    \frac{1}{\sum_{j=1}^n\beta^j}\le\alpha\le
    \frac{-\beta^{n+1}+\beta^n+2\beta-1}{\sum_{j=1}^n\beta^j},
    \]
    or (symmetrical case)
    \[
    2-\beta+\frac{\beta^{n+1}-\beta^n-2\beta+1}{\sum_{j=1}^n\beta^j}
    \le\alpha\le
    2-\beta-\frac{1}{\sum_{j=1}^n\beta^j},
    \]
    then $T$ is not topologically transitive. The sets described by the above inequalities are visible in Fig.~\ref{fig:triangleTrans} as white lens-shaped areas whose upper vertices lie on red curves. However, the figure suggests that, in addition to the areas described by the above inequalities, there are other areas (also lens-shaped) distributed fractally (as in Cantor set constructions) between them. For example, observe the part of the triangle between the largest and the second-largest white lens.

    \item If $\beta\ge\sqrt{2}$, then $T$ is topologically transitive. This corresponds to the smaller blue triangle above the horizontal red line in Fig.~\ref{fig:triangleTrans}. In fact, the authors of~\cite{mixing} proved more than transitivity here, namely the mixing property (with a one-point exception).

    \item If $\sqrt{2}>\beta\ge\sqrt[3]{2}$ and $1/(\beta^2+\beta)>\alpha\ge0$, or           $2-\beta\ge\alpha>2-\beta-1/(\beta^2+\beta)$, then $T$ is topologically transitive. This corresponds to the blue area in the strip of the triangle $\mathcal{T}$ between the lines $\beta=\sqrt{2}$ and $\beta=\sqrt[3]{2}$ (excluding, of course, the part of the largest white lens). Here, we also have mixing (with two points exception).
\end{enumerate}

More interestingly, Fig.~\ref{fig:triangleTrans} allows us to formulate a new conjecture concerning $\beta$-transformations, which has not yet been proved.
\begin{conj}
For all parameters from the triangle 
$\mathcal{T}$ satisfying the conditions
\[
    \big(\beta\le\sqrt{2}
    \quad\text{and}\quad
    \alpha\le1-1/\beta\big)
\]
or
\[
    \big(\beta\le\sqrt{2}
    \quad\text{and}\quad
    \alpha\ge1+1/\beta-\beta\big),
\]
i.e., belonging to the union of 
the two blue curvilinear triangles
restricted by the red curves,
the corresponding $\beta$-transformations
are topologically transitive.
\end{conj}

\begin{rem}
In fact, we believe that 
$\beta$-transformations
are topologically transitive for all parameters
outside the white lenticular regions
from Fig.~\ref{fig:triangleTrans},
but unfortunately at the moment we cannot provide
the precise formal description 
of this whole white area. 
\end{rem}
\begin{figure*}[!htb]
    \centering
    \includegraphics[scale=0.44, trim=0mm 5mm 5mm 0mm]{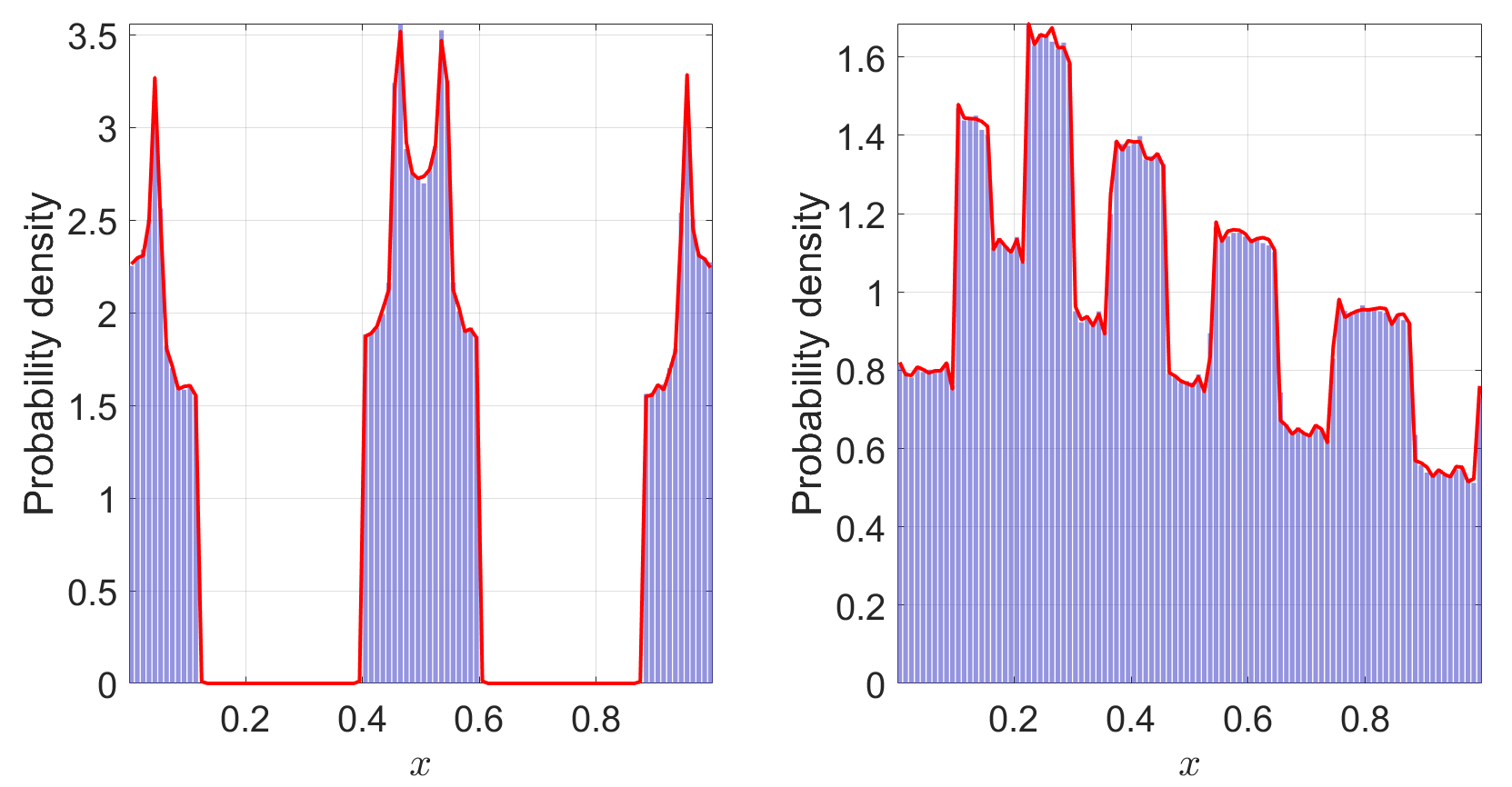}
    \caption{Probability density functions and histograms
    for nontransitive (left panel) and transitive 
    (right panel) $\beta$-transformations. Parameter
    values: $\alpha=0.4$, $\beta=1.2$ (left) and
    $\alpha=0.1$, $\beta=1.2$ (right).}
    \label{fig:probdensity}
\end{figure*}
Finally, let us explain briefly why our numerical transitivity 
test works. It seems that it is related to 
the existence of an absolutely continuous invariant probability 
measure (acip for short). Namely, if a map on the unit interval
admits an ergodic absolutely continuous measure, then every
computer-generated orbit consists of contiguous segments, 
where each segment is close to a segment on a dense real orbit 
(for more detailed and thorough explanations of 
the impact of acip measures on computer orbits 
consult \cite{boyarsky1986,friedman1988,gora1988}
and references therein).
The key fact here is the following important
result, which is true for a wide range of expanding Lorenz maps
(see the proof of Theorem 5.5 in~\cite{cnv_model1}).

\begin{proposition}
\label{prop:acip}
Every piecewise $C^2$ expanding Lorenz map 
has a unique ergodic acip 
measure, and the support of this measure 
is a finite union of closed intervals.
\end{proposition}

Now Proposition~\ref{prop:acip} leads 
to the following alternative,
which can be observed in Fig.~\ref{fig:probdensity}
showing the probability density function of acip measures.
Either the support of the acip measure is
the whole domain (one full interval)
as on the right panel of Fig.~\ref{fig:probdensity},
or it is a proper subset of the domain,
which does not exhaust the entire domain
(a finite union of disjoint intervals),
as on the left panel of Fig.~\ref{fig:probdensity}.
In the first case an expanding Lorenz map
is topologically transitive (in the domain), 
and in the second one
it is not.  
\section{Numerical LEO}\label{sec:numleo}

This section is devoted to the study of \emph{numerical LEO},
which is the computational equivalent of 
\emph{topological locally eventually onto} property 
(see the formal definition below).
It is well known that the topological LEO property implies
topological transitivity, but not vice versa.
Namely, the $\beta$-transformation 
$f(x)=\sqrt{2}x+\frac{2-\sqrt{2}}{2} \pmod 1$
from $[0,1)$ into itself is topologically transitive,
but not topologically LEO (for more details see~\cite{mixing}).
Roughly speaking, topological transitivity represents 
a \emph{weaker} version of chaotic behavior, 
while topological LEO represents a \emph{stronger} version.
We believe that this correspondence also occurs 
at the numerical level, i.e., that the numerical LEO property
should be stronger than the numerical transitivity,
assuming that both algorithms work correctly.

Recall the definition of the topological LEO property.
Assume $f\colon X\to X$ is from a metric space into itself.
Then $f$ is called 
\emph{topologically locally eventually onto} 
(\emph{topologically LEO} for short) 
if for every nonempty open set 
$U$ there is $n\in\mathbb{N}$ such that $f^n(U)=X$.
If $f$ is not LEO we will call it
\emph{topologically nonLEO}.
Note that, by definition, topological LEO implies 
topological mixing.

\vspace{1mm}
\textbf{Algorithm of numerical LEO.}
Below we present the main idea of the numerical LEO test,
which is completely and directly based on 
the LEO property definition.

\begin{enumerate}
    \item The algorithm tests whether the iterates 
    of a Lorenz-like map $f$ eventually 
	cover the entire interval \(I=[0,1)\), 
    even when starting from small subintervals of \(I\).
    \item The main procedure divides the interval 
    \(I\) into smaller parts 
	(subintervals) and verifies if each subinterval, 
	under repeated application of the map, expands 
    to fill the whole domain.
    \item For each subinterval, the algorithm tracks 
    how it transforms through many iterations, 
    considering the effect of a discontinuity 
    in the Lorenz-like map.
    \item The evolution of each interval is approximated 
    by computing images of its endpoints and adjusting for cases 
    where the image contains the discontinuity point.
    \item After each iteration, overlapping image 
    intervals are merged, and the process stops early 
    if full coverage is achieved.
    \item If all initial subintervals succeed in covering 
    the entire domain, the test confirms that the map 
    exhibits the numerical LEO property.
\end{enumerate}

As we will see, implementing these guidelines 
in pseudocode leads to a much more complex algorithm 
than in the case of numerical transitivity. 
Due to its much larger size, we have divided 
this big algorithm into four submodules. Namely,
Algorithm~\ref{alg:merge} merges overlapping
intervals (by default, images of a Lorenz map).
Algorithm~\ref{alg:image} computes the image of 
the fixed subinterval under some fixed iteration
of a Lorenz map.
Algorithm~\ref{alg:cover} checks 
if a fixed small subinterval covers the entire
domain after some fixed number of iterations 
of a Lorenz map.
Finally, the input of the main module, i.e., 
Algorithm~\ref{alg:LeoTest}, is a Lorenz map $f$ 
and three numbers $c$ (discontinuity point) and $a$, $b$ 
(endpoints of the domain of $f$), and
the output is \texttt{true} (i.e., $1$) if $f$ is numerically LEO 
and \texttt{false} (i.e., $0$) if $f$ is not.

\begin{algorithm}[!htb]
{\small
\caption{\textsc{Merging Finite List of 
Overlapping Intervals}}\label{alg:merge}
\begin{algorithmic}[1]
\Function{MergeIntervals}{$\text{\it{intervals}}$}
    \Comment{Merges overlapping intervals into a consolidated list}
    
    \State \textsc{Sort} $\text{\it{intervals}}$ by the first value in each interval 
    \Comment{Ensure intervals are processed in ascending order}
    
    \State $\text{\it{merged}} \gets [\text{\it{intervals}}[1]]$ 
    \Comment{Initialize merged list with the first interval}
    
    \For{$i = 2$ to length of $\text{\it{intervals}}$} 
        \Comment{Iterate through remaining intervals}
        
        \State $\text{\it{current}} \gets \text{\it{intervals}}[i]$ 
        \Comment{Current interval to process}
        
        \State $\text{\it{last}} \gets \text{\it{merged}}[\text{end]}$ 
        \Comment{Last interval in the merged list}
        
        \If{$\text{\it{current}}[1] \leq \text{\it{last}}[2]$} 
            \Comment{Check if intervals overlap}
            \State $\text{\it{merged}}[\text{end}][2] \gets \max(\text{\it{last}}[2], \text{\it{current}}[2])$ 
            \Comment{Merge intervals by updating the endpoint}
        \Else
            \Comment{No overlap, add current interval to merged list}
            \State \textsc{Append} $\text{\it{current}}$ to $\text{\it{merged}}$
        \EndIf
    \EndFor
    
    \State \Return $\text{\it{merged}}$ 
    \Comment{Return the merged list of intervals}
\EndFunction
\end{algorithmic}}
\end{algorithm}

\begin{algorithm}[!htb]
{\small
\caption{\textsc{Image of Interval under Iterated Map as
Union of Intervals}}\label{alg:image}
\begin{algorithmic}[1]
\Function{Image}{$T$, $c$, $x$, $y$, $a$, $b$}
    \Comment{Computes the image of interval $[x, y]$ under the $l$-th iteration of a Lorenz map 
    $T\colon[a,b]\to[a,b]$ with a discontinuity at $c$.}
    \State $l \gets$ number of iterations
    \Comment{Fix some big number of iterations}
    \State $\text{\it{intervals}} \gets [[x, y]]$ 
    \Comment{Initialize with the given interval}
    
    \For{$i = 1$ to $l$} 
        \State $\text{\it{new\_intervals}} \gets []$
        \Comment{Initialize with the empty list}
        \For{each $\text{\it{interval}}$ in $\text{\it{intervals}}$}
            \State $x1 \gets \text{\it{interval}}[1]$, $x2 \gets \text{\it{interval}}[2]$
            \State $T1 \gets T(x1)$, $T2 \gets T(x2)$
            \Comment{Compute the image of endpoints}
            
            \If{$x1 < c$ and $c < x2$} 
                \Comment{Interval crosses discontinuity}
                \State $\text{\it{new\_intervals}} \gets
                [\text{\it{new\_intervals}}; [a, T2]; [T1, b]]$
            \ElsIf{$x2 = c$} 
                \Comment{Interval ends at discontinuity}
                \State $\text{\it{new\_intervals}} \gets
                [\text{\it{new\_intervals}}; [T1, b]]$
            \Else
                \Comment{Interval does not cross discontinuity}
                \State $\text{\it{new\_intervals}} \gets
                [\text{\it{new\_intervals}}; [T1, T2]]$
            \EndIf
        \EndFor
        \State $\text{\it{intervals}} \gets \textsc{MergeIntervals}(\text{\it{new\_intervals}})$ 
        \Comment{Merge overlapping intervals}
        
        \If{there is only one interval $[a,b]$}
            \Comment{Break if the entire space is covered}
            \State \textbf{break}
        \EndIf
    \EndFor
    
    \State \Return $\text{\it{intervals}}$
\EndFunction
\end{algorithmic}}
\end{algorithm}

\begin{algorithm}[!htb]
{\small
\caption{\textsc{Numerical Covering Test for Domain Coverage}}
\label{alg:cover}
\begin{algorithmic}[1]
\Function{NumCoverTest}{$T$, $c$, $x$, $y$,
$a$, $b$}
    \Comment{Checks if the interval $[x, y]$ covers the entire domain $[a,b]$ after some big number of iterations of a Lorenz map  
    $T\colon[a,b]\to[a,b]$.}
    \State $\text{\it{logValue}} \gets$ \textsc{false}
    \Comment{Initialize logical flag}
    
    \State $\text{\it{intervals}} \gets \textsc{Image}(T, c, x, y, a, b)$ 
    \Comment{Compute the image of the interval $[x, y]$
    after some big number of iterations}
    
    \If{there is only one interval $[a,b]$ in $\text{\it{intervals}}$} 
        \Comment{Check if the entire space is covered}
        \State $\text{\it{logValue}} \gets$ \textsc{true}
    \EndIf
    
    \State \Return $\text{\it{logValue}}$
\EndFunction
\end{algorithmic}}
\end{algorithm}

\begin{algorithm}[!htb]
{\small
\caption{\textsc{Numerical Leo Test for Domain Coverage}}
\label{alg:LeoTest}
\begin{algorithmic}[1]
\Function{NumLeoTest}{$T$, $c$, $a$, $b$}
    \Comment{Checks if each small subinterval from
    partition covers the entire space $[a,b]$ after some big number of iterations of a Lorenz map $T$
    with a discontinuity at $c$.}
    \State $\text{\it{logValue}} \gets$ \textsc{false}
    \Comment{Initialize logical flag}
    \State $m \gets$ number of subdivisions
    \Comment{Fix a number of parts of division}
    
    \State $dm \gets (b-a)/m$ 
    \Comment{Divide the interval $[a,b]$ into $m$ equal subintervals}
    
    \For{$i = 1$ to $m$} 
        \State $x \gets a+(i-1) \cdot dm$
        \Comment{Beginning of subinterval}
        \State $y \gets a+i \cdot dm$
        \Comment{End of subinterval}
        \State $\text{\it{logTest}} \gets \textsc{NumCoverTest}(T, c, x, y, a, b)$
        \Comment{Test subinterval $[x,y]$}
        
        \If{$\text{\it{logTest}} = 0$} 
            \Comment{If any subinterval fails, return false}
            \State \Return \textsc{false}
        \EndIf
    \EndFor
    
    \State $\text{\it{logValue}} \gets$ \textsc{true}
    \Comment{If all subintervals pass, return true}
    \State \Return $\text{\it{logValue}}$
\EndFunction
\end{algorithmic}}
\end{algorithm}

\section{Comparison of numerical transitivity 
and numerical LEO tests}
\label{sec:comparison}

This section is devoted to the comparison
of numerical transitivity and numerical LEO algorithms
and tests. The \emph{numerical transitivity test} and 
the \emph{numerical LEO test} are two numerical algorithms 
designed to analyze different but related dynamical properties 
of maps on the interval $[a, b)$. 
While both aim to characterize complex behavior 
such as mixing or chaos, they differ significantly 
in their assumptions, approach, and interpretation.

The \emph{numerical transitivity test} simulates 
the trajectory (or time series) of a single point 
under iteration of a map $f\colon [a, b) \to [a, b)$. 
It checks whether the orbit of a randomly chosen initial point
becomes dense in the domain by verifying that 
it visits all bins in a fine partition of the interval. 
This test is probabilistic in nature: 
it uses multiple random trials and a fixed number of iterations, 
discarding a transient portion of the trajectory. 
If at least one trial results in a trajectory 
that visits all subintervals, the test is considered passed.

In contrast, the \emph{numerical LEO test} is a deterministic 
and more structured approach tailored for Lorenz-type maps 
with a known discontinuity at a point $c$. 
It examines whether every small subinterval 
of the domain eventually maps to cover the entire interval 
$[a, b]$ after a fixed number of iterations. 
This is achieved by computing successive images 
of each subinterval and merging overlapping intervals. 
The test only passes if \emph{all} subintervals 
eventually cover the domain, indicating 
a strong form of topological mixing.

Overall, intendedly, the numerical transitivity 
test provides a quick, 
sample-based indication of chaotic behavior, 
while the numerical LEO test offers a more rigorous 
and exhaustive check of domain covering properties. 
The former is suitable for general maps, whereas 
the latter assumes a specific structure and 
seems to be more robust to missing rare dynamical features.
Moreover, the numerical transitivity test
is many times faster, i.e., it requires less computational 
time to execute, is more efficient in terms 
of algorithmic complexity and, finally,
uses fewer resources (CPU, memory).

\begin{figure*}[htbp]
    \centering
    \includegraphics[scale=0.24,trim=0mm 0mm 0mm 0mm]{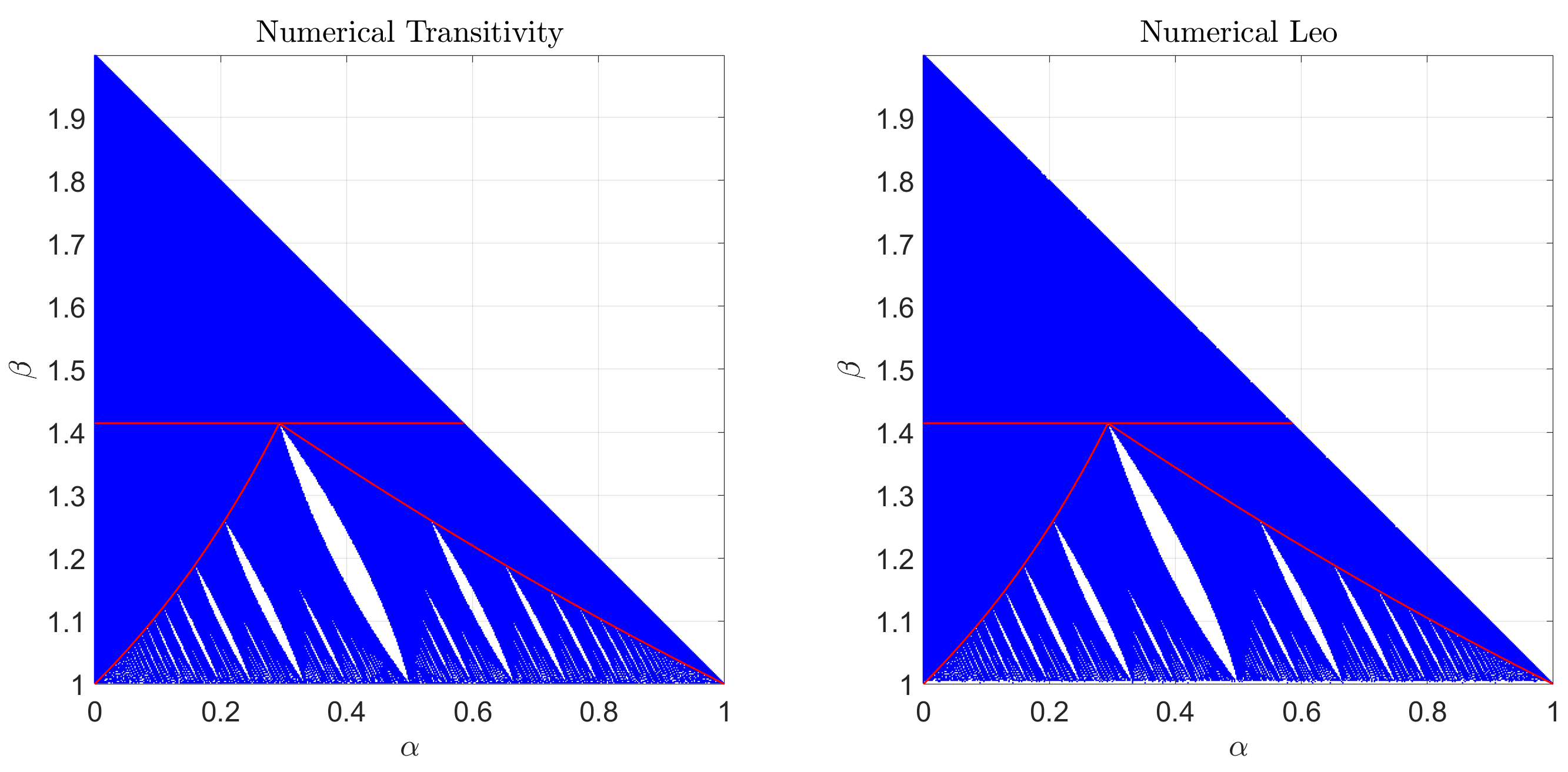}
    \caption{
    Comparison of results of numerical transitivity 
    and numerical LEO tests for classical $\beta$-transformations
    in the $\alpha$–$\beta$ parameter plane.
    Meshgrid = 500.
    }
    \label{fig:numandtrans}
\end{figure*}
And now, to a big surprise, both tests give extremely similar results in the $\beta$-transformation class.
Namely, Figure~\ref{fig:numandtrans} presents a comparison of the graphical representation of the results of both tests for selected, evenly distributed values (meshgrid $500$) from the triangle $\mathcal{T}$ of $\alpha$-$\beta$ parameters for $\beta$-transformations.
It seems that, with the exception of a narrow strip just above the $\alpha$-axis, where the numerical LEO algorithm does not work (due to the slope of the line close to $\beta=1$, the number of iterations needed to cover the entire domain increases rapidly and quickly exceeds any fixed value stored in the algorithm), the results of both tests are almost identical.

To be honest, the blue color of the triangle $\mathcal{T}$ above the horizontal red line $\beta=\sqrt{2}$ is no surprise, but rather a confirmation of the consistency of the tests with the theoretical results.
Namely, the following result proved in~\cite{allbutone} guarantees topological LEO and, in consequence, also topological transitivity for parameters from this area.

\begin{thm}\label{thm:square}
Assume $f$ is an expanding Lorenz map and
\[
\beta_f = \inf{\{f'(x)\mid x\in[0,1)\setminus F\}},
\]where $F$ is the set of non-differentiability points.
Let
\[
f_0(x) = \sqrt{2}x + \frac{2-\sqrt{2}}{2} \pmod 1.
\]
If $\beta_f \ge \sqrt{2}$ and $f\neq f_0$, then $f$ is topologically LEO.
\end{thm}However, it should be emphasized that apart from the above, no other rigorous results concerning the topological LEO property for $\beta$-transformations are currently known. Only some partial results concerning the topological mixing property have been discovered (for more details see~\cite{mixing}).
In this way, the similarity between the left and right panels of Figure~\ref{fig:numandtrans} in the area below the red line $\beta=\sqrt{2}$ at first sight seems puzzling.

What do our numerical simulations presented in Figure~\ref{fig:numandtrans} suggest? They suggest that for $\beta$-transformations \textsc{topological transitivity} = \textsc{topological LEO}.
In other words, the following alternative holds: either a $\beta$-transformation is topologically LEO (blue area) or it is topologically nontransitive (white area). Note that, in general, this is \emph{not true} from a formal point of view, because there are $\beta$-transformations that are topologically transitive and not topologically LEO (for example, the map $f_0$ from Theorem~\ref{thm:square}).
However, it seems that such a situation (topological transitivity and topological nonLEO) is actually very rare and exceptional and therefore it is not visible in our drawing.
With this in mind, Figure~\ref{fig:numandtrans} suggests the following hypothesis not yet proven.

\begin{conj}
For all $\beta$-transformation parameter values from the triangle $\mathcal{T}$ except for countably many, we have \textsc{topological transitivity} = \textsc{topological LEO}, i.e., a $\beta$-transformation is topologically transitive if and only if it is topologically LEO.
In consequence, this equality holds with probability $1$, i.e., the subset of parameters from the triangle $\mathcal{T}$ for which this equality does not hold has Lebesgue measure $0$.
\end{conj}
\begin{rem}
Observe that for now we cannot provide a complete and precise description of the set of parameters for which the $\beta$-transformation is topologically transitive, nor the set for which it is topologically LEO.
Nevertheless, our hypothesis makes sense, because it does not refer to these sets in any way. Let us also mention another interesting hypothesis related to the topological LEO property: if an expanding Lorenz map is topologically mixing then it is also topologically LEO.
\end{rem}

\begin{figure}[!htb]
    \centering
    \includegraphics[scale=0.19, trim=0mm 0mm 0mm 0mm]{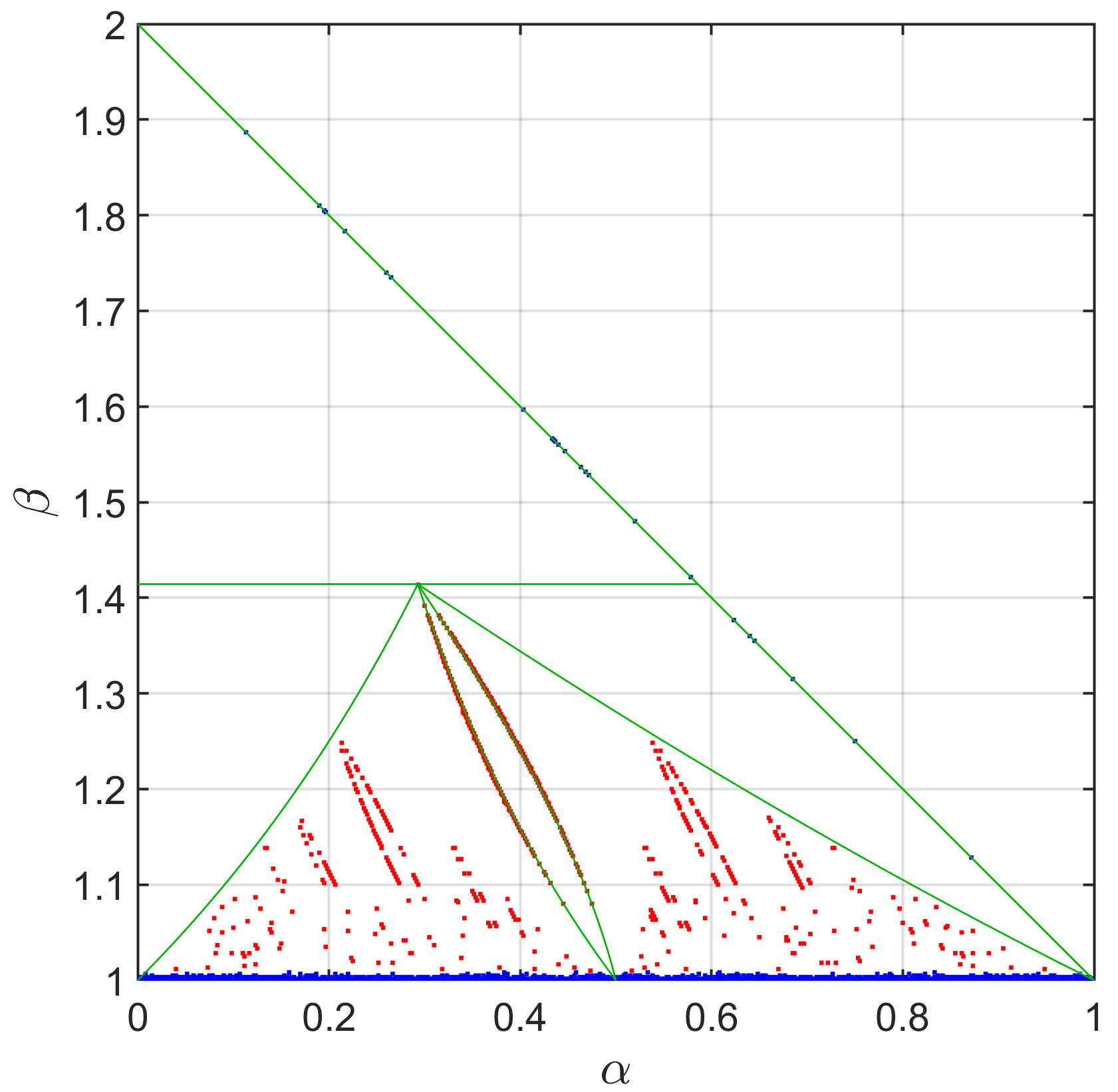}
    \caption{
    Differences between numerical transitivity 
    and numerical LEO tests results for classical 
    $\beta$-transformations in the $\alpha$-$\beta$ parameter plane.
    Blue points: numerical transitivity and not numerical LEO,
    Red points: numerical LEO and not numerical transitivity.
    Meshgrid = 600.
    }
    \label{fig:difftriangle}
\end{figure}
Although, as shown in Figure~\ref{fig:numandtrans}, the results of the numerical transitivity and numerical LEO tests are very similar, there are points in the triangle $\mathcal{T}$ where the results of both tests do not match.
Let us take a closer look at these points, which are visible in Figure~\ref{fig:difftriangle}.
In short, it seems that what we observe are mainly \emph{numerical artifacts} that arise due to the limitations of our numerical methods.
The two main reasons for these limitations are the consideration of only a finite division of the domain into subintervals (in the case of numerical transitivity) and a finite number of iterations (in the case of numerical LEO).

The red dots in Figure~\ref{fig:difftriangle} correspond to the parameter values $(\alpha,\beta)$ for which the $\beta$-transformation is not numerically transitive but is numerically LEO.
Note that formally there are no maps that are topologically LEO and topologically nontransitive (see the beginning of Section~\ref{sec:numleo}).
Moreover, the red dots appear on the edges of white lenses from Figure~\ref{fig:triangleTrans}.
Since, from a formal point of view, the points on the edges of the lenses correspond to maps that are not topologically transitive (see Section~\ref{subsec:transsim}), probably due to computational limitations, the numerical LEO test gives here a ``wrong'' result that is inconsistent with theoretical predictions.

In turn, the blue dots correspond to the parameter values for which the $\beta$-transformation is numerically transitive but not numerically LEO.
Probably, once again, in most cases, the numerical LEO test ``fails'' here, i.e. it is inconsistent with the theoretical LEO.
In particular, the results of the numerical LEO test are not very reliable for points located very close to the base of the triangle $\mathcal{T}$, because at these points the slope of the plot is close to $1$ and, in consequence, a very large number of iterations is needed for the small initial interval to cover the entire domain.
Similarly, the blue dots on the green main diagonal seem to be numerical artifacts.

\section{Courbage--Nekorkin--Vdovin model}
\label{sec:cnv}

The Courbage--Nekorkin--Vdovin (CNV) model was first developed in 2007
as a two-dimensional discrete system to describe 
how neurons produce repeated electrical spikes,
which are known as action potentials~\cite{courbage2007,courbage2010}. 
Both regular and irregular (chaotic) spike patterns are seen 
in this 2D representation, just like in real neurons. 
To make it easier to analyze and simulate, 
both the authors of~\cite{courbage2007} and subsequent researchers 
(e.g.,~\cite{cnv_model}) developed a one-dimensional version of the model 
that preserves the essential features that lead to these complex behaviors.
In this 1D model, the only variable is potential, 
while the second variable of the 2D system becomes 
the main parameter, the so-called recovery parameter.
This more straightforward 1D model is better suited 
for mathematical research while still capturing crucial neuronal dynamics 
like spiking and chaos. We use both piecewise linear and nonlinear versions 
of this 1D CNV model, where the linear version allows for more direct
theoretical analysis, while the nonlinear version, which is based on 
a cubic-like map, exhibits even richer dynamics and is suitable 
for studying properties such as topological transitivity, topological mixing, 
and topological LEO through numerical simulations~\cite{cnv_model1,multilevel}.

\subsection{Description of CNV model}

The one-dimensional \emph{Courbage--Nekorkin--Vdovin} 
(\emph{CNV}, for short)
model of a neuron is described by the formula
\begin{equation}
\label{eq:mapg}
x_{n+1}=g(x_n)
=x_n+F(x_n)-\alpha-\beta H(x_n-d),
\end{equation}
where
\begin{itemize}
  \item $x$ is the \emph{membrane potential} of the neuron
        (the main model variable),
  \item $H(x)$ is the Heaviside step function, i.e.,
        \[
        H(x)=
        \begin{cases} 
        1, & \text{if } x \geq 0,\\
        0, & \text{if } x < 0,
        \end{cases}
        \]
  \item $F(x)$ is a specially chosen continuous function 
        whose plot has the shape of an upside-down, reversed
        $\boldsymbol{N}$ letter,
  \item $\alpha$, $\beta$, and $d$ are model parameters.
\end{itemize}

Note that, in general, we have four parameters:
\begin{enumerate}
    \item $\alpha$ — the main parameter
    of the model, called the \emph{recovery} parameter,
    \item $\beta$ — the \emph{jump} of the scaled 
    step function and also 
    the length of the \emph{invariant
    interval} (if it exists),
    \item $d$ — the value of the \emph{discontinuity}
    point,
    \item the function $F$, which may depend on
    additional parameters.
\end{enumerate}

\begin{figure*}[!htb]
    \centering
    \includegraphics[scale=0.31,clip=true,
     trim=0mm 0mm 0mm 0mm]{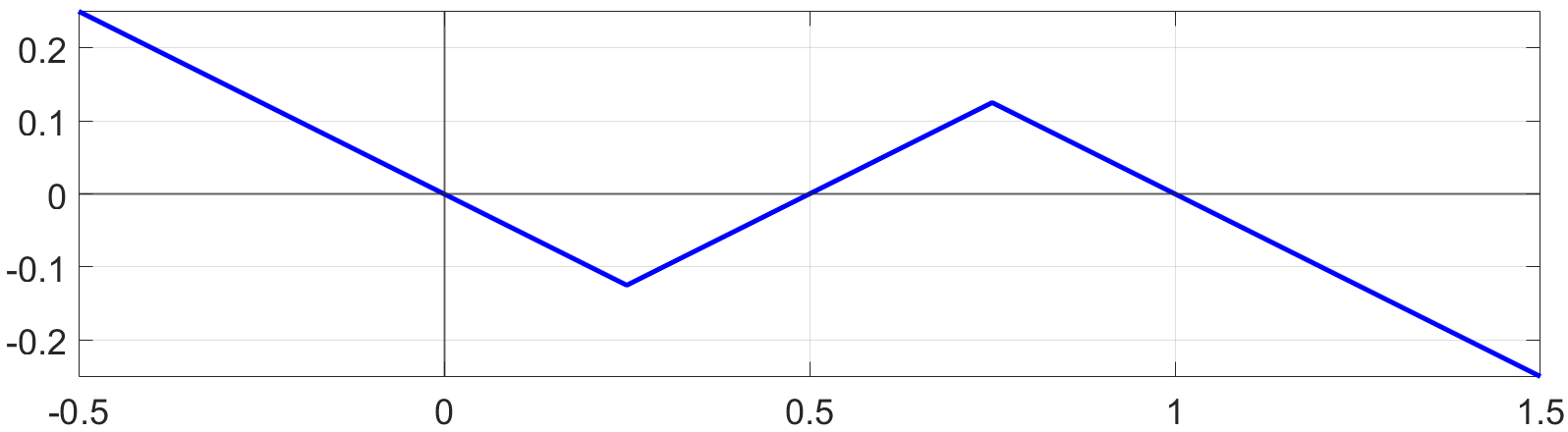}
    \caption{Plot of $F$ function
    in the plCNV case.}
    \label{fig:model}
\end{figure*}

We usually consider two main versions of the 1D CNV model:
\begin{itemize}
    \item a \emph{piecewise linear} case (plCNV, for short), 
    when $F(x)$ is a piecewise linear continuous
    function defined as follows:
    \[
    F(x) =
    \begin{cases} 
    -m_0 x, & \text{if } x \leq J_{\min},\\[3pt]
    m_1(x-a), & \text{if } J_{\min} \leq x \leq J_{\max},\\[3pt]
    -m_0(x-1), & \text{if } x \geq J_{\max},
    \end{cases}
    \]
    where 
    \[
    m_0, m_1>0,\quad
    1 > a > 0,\]
    \[
    J_{\min}=\frac{a m_1}{m_0+m_1},\quad
    J_{\max}=\frac{m_0+a m_1}{m_0+m_1},
    \]
    \item a \emph{nonlinear} case (nlCNV, for short),
    when $F(x)=\mu x(x-a)(1-x)$
    with $0<a<1$ and $\mu>0$.
\end{itemize}

Observe that in both cases the function $F$
is continuous and piecewise monotonic, 
but the function $g$ is discontinuous at $d$. 

\begin{figure*}[!htb]
    \centering
    \includegraphics[scale=0.33,clip=true,trim=2mm 0mm 0mm 0mm]{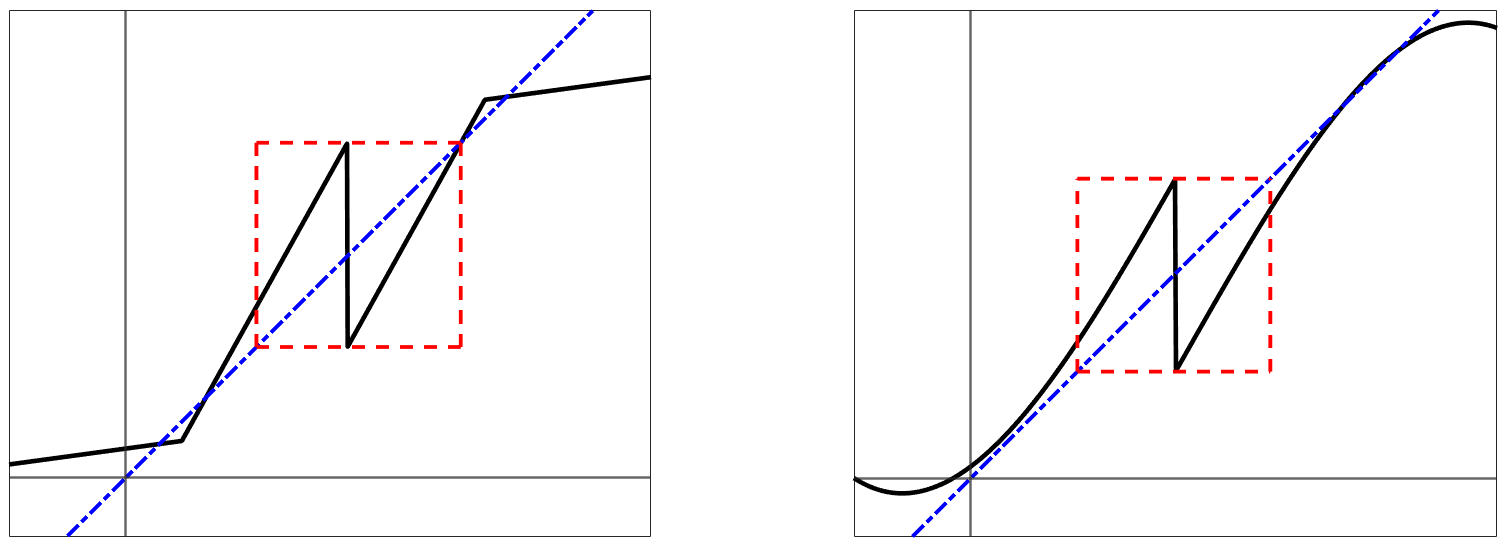}
    \caption{Example plots of plCNV (left)
    and nlCNV (right) model functions.
    The red boxes indicate invariant intervals.}
    \label{fig:model2}
\end{figure*}

\subsection{Conditions for the existence
of the invariant interval}

Now we examine the existence of
the invariant interval, 
on which the 1D plCNV function
is a $\beta$-transformation
(resp. the 1D nlCNV function
is an expanding Lorenz map).
Such an invariant interval is indicated 
in Fig.~\ref{fig:model}. 
The endpoints of this interval are given,
in both cases, 
by the formulas
\begin{align*}
 b &= \lim_{x\to d^+}g(x)=g(d)=d+F(d)-\alpha-\beta,\\
 c &= \lim_{x\to d^-}g(x)=d+F(d)-\alpha.
\end{align*}
Note that using a linear 
transformation of variables
(see above formulas),
instead of $\alpha$ and $\beta$
we can treat $b$ and $c$ as the main
parameters of the model (with all other
parameters fixed). Below we provide
conditions for the existence of
the invariant interval in terms
of $b$ and $c$.
\paragraph{plCNV case}

In Table~\ref{table:plsixcond} we present 
necessary and sufficient conditions for the existence
of the invariant interval on which
a plCNV model map is a $\beta$-transformation,
and we explain their geometric meaning.

\begin{table}[!htb]
    \centering
    \caption{Conditions for the existence
    of the invariant interval (plCNV).}
    \label{table:plsixcond}
    \begin{tabular}{@{}lll@{}}
    \toprule
    No. & Condition & Geometric meaning\\
    \midrule
    1. & $J_{\min}\le b$ & Monotonicity of the left branch\\
    2. & $c\le J_{\max}$ & Monotonicity of the right branch\\
    3. & $b<d$ & Proper position of discontinuity\\
    4. & $d<c$ & Proper position of discontinuity\\
    5. & $g(b)\ge b$ & Map from $[b,c)$ into itself\\
    6. & $g(c)< c$ & Map from $[b,c)$ into itself\\
    \bottomrule
    \end{tabular}
\end{table}

\paragraph{nlCNV case}

In turn, Table~\ref{table:nlsixcond} provides 
similar conditions for the existence
of the invariant interval in the case
of the nlCNV model.
Let $x_{\min}$ and $x_{\max}$ denote 
the points of local
minimum and maximum of the cubic polynomial $F$.
An immediate calculation gives
\[
x_{\min}=
\frac{a+1-\sqrt{a^2-a+1}}{3},
\qquad
x_{\max}=
\frac{a+1+\sqrt{a^2-a+1}}{3}.
\]

\begin{table}[!htb]
    \centering
    \caption{Conditions for the existence
    of an invariant interval (nlCNV).}
    \label{table:nlsixcond}
    \begin{tabular}{@{}lll@{}}
    \toprule
    No. & Condition & Parametric form\\
    \midrule
    1. & $x_{\min}<b$ & $\beta<d+F(d)-x_{\min}-\alpha$\\
    2. & $c<x_{\max}$ & $\alpha>d+F(d)-x_{\max}$\\
    3. & $b<d$ & $\beta>F(d)-\alpha$\\
    4. & $d<c$ & $\alpha<F(d)$\\
    5. & $g(b)\ge b$ & $\alpha\le F(d+F(d)-\alpha-\beta)$\\
    6. & $g(c)< c$ & $\beta> F(d+F(d)-\alpha)-\alpha$\\
    \bottomrule
    \end{tabular}
\end{table}

In what follows, we will use the above conditions 
in all our numerical simulations to determine 
a set of parameters (usually as a region 
in the $\alpha$--$\beta$ or $b$--$c$ parameter plane) 
for which we will perform numerical transitivity 
and LEO tests. This set of parameters will be
represented in all figures as a colored area.
Finally, note that since $[b,c)\subset(x_{\min},x_{\max})$,
we get $F'(x)>0$ and, in consequence, $g'(x) > 1$ 
for $x\in[b,c)$, which is needed in the definition 
of an expanding Lorenz map.

\subsection{Summary of previous results on the CNV model}

Let us recall that both models—the piecewise linear (1D plCNV) 
and nonlinear (1D nlCNV)—have been studied
in detail in recent papers~\cite{cnv_model,cnv_model1,multilevel}. 
In particular, in~\cite{cnv_model} the authors 
analyze existence, position and stability 
of fixed points, giving explicit conditions 
with respect to model parameters. 
After restricting the model 
to the invariant interval (beyond this interval 
the dynamics is trivial), they establish the
conditions for Devaney chaos and describe metric properties 
of this model, among other things, existence and form 
of the absolutely continuous invariant probability measure. 
In addition, the itineraries of periodic orbits
are linked with patterns of 
spike trains fired by the model map.
Note that some of the obtained results are 
due to the fact that the 1D plCNV model 
map is a $\beta$-transformation, 
not a general expanding Lorenz map. 

In turn, in~\cite{cnv_model1} the authors 
transfer some of the results mentioned above
to the nlCNV model and identify differences 
between the piecewise linear and nonlinear cases.
Precisely, the article presents a detailed analysis of 
the 1D nlCNV model, identifying parameter regions 
where the system exhibits expanding Lorenz map behavior. 
This framework enables the application of Lorenz map theory 
to establish sufficient conditions for chaos, analyze rotation intervals, 
and characterize periodic orbit itineraries, revealing 
the quite intricate structure of spike patterns. 
Moreover, numerical simulations support 
these theoretical results, providing interesting examples 
of rotation intervals and orbit itineraries for selected parameters. 
Finally, the authors of~\cite{cnv_model1} extend 
the discussion to the non-autonomous discrete case, 
presenting voltage time series for time-varying inputs.

In the later work~\cite{multilevel}, the authors develop
the rotation interval analysis to study 
the regularity of orbits in expanding Lorenz maps 
and link these properties to the CNV model. 
They point out that certain features 
of CNV model maps play a key role 
in influencing long-term firing patterns of neurons.
Namely, they examine the structure and characteristics 
of periodic orbit itineraries for expanding Lorenz maps. 
Specifically, they demonstrate that the periodic orbits 
in these maps form two distinct cascades (Stern--Brocot 
and Geller--Misiurewicz) that are intimately related to 
the Farey tree of rational rotation numbers 
associated with the map. 
These results are applied to the 1D CNV neuron model, 
providing a detailed description of the multilevel regularity 
in its periodic spiking patterns.
\section{Numerical transitivity and numerical LEO in the CNV model}
\label{sec:transleocnv}

In this section we present the results of numerical 
simulations concerning the occurrence of 
the numerical transitivity (the first subsection)
and the numerical LEO (the second subsection)
in the 1D CNV neuron model.
Each subsection is divided into two parts
related to two versions of the model:
plCNV and nlCNV. Recall that
in both cases, the model map has the form
\[
g(x)
= x + F(x) - \alpha - \beta H(x - d),
\]
where $\alpha$ and $\beta$ are
the main parameters of the model.
Our main goal here is to analyze how our numerical properties 
(transitivity and LEO) depend on these parameters.
However, it turns out that we obtain more graphically 
readable results if we replace the $\alpha$ and $\beta$ 
parameters with the $b$ and $c$ parameters
corresponding to the endpoints of the invariant interval
for $g$. As we explained earlier
in Section~\ref{sec:cnv}, this parameter conversion
is possible thanks to the use of the linear change of 
coordinates given by the formulas for $b$ and $c$.
Moreover, additional symmetry with respect to one 
of the diagonals visible in all figures in this section
is related to the choice of the discontinuity
point $d$ as $(x_{\min} + x_{\max}) / 2$.
We believe that the presented simulations offer 
an essential insight into how different levels
of chaotic behavior emerge
and affect neuronal responses.

\begin{figure*}[!htb]
    \centering
    \includegraphics[scale=0.21, trim=5mm 5mm 0mm 0mm]{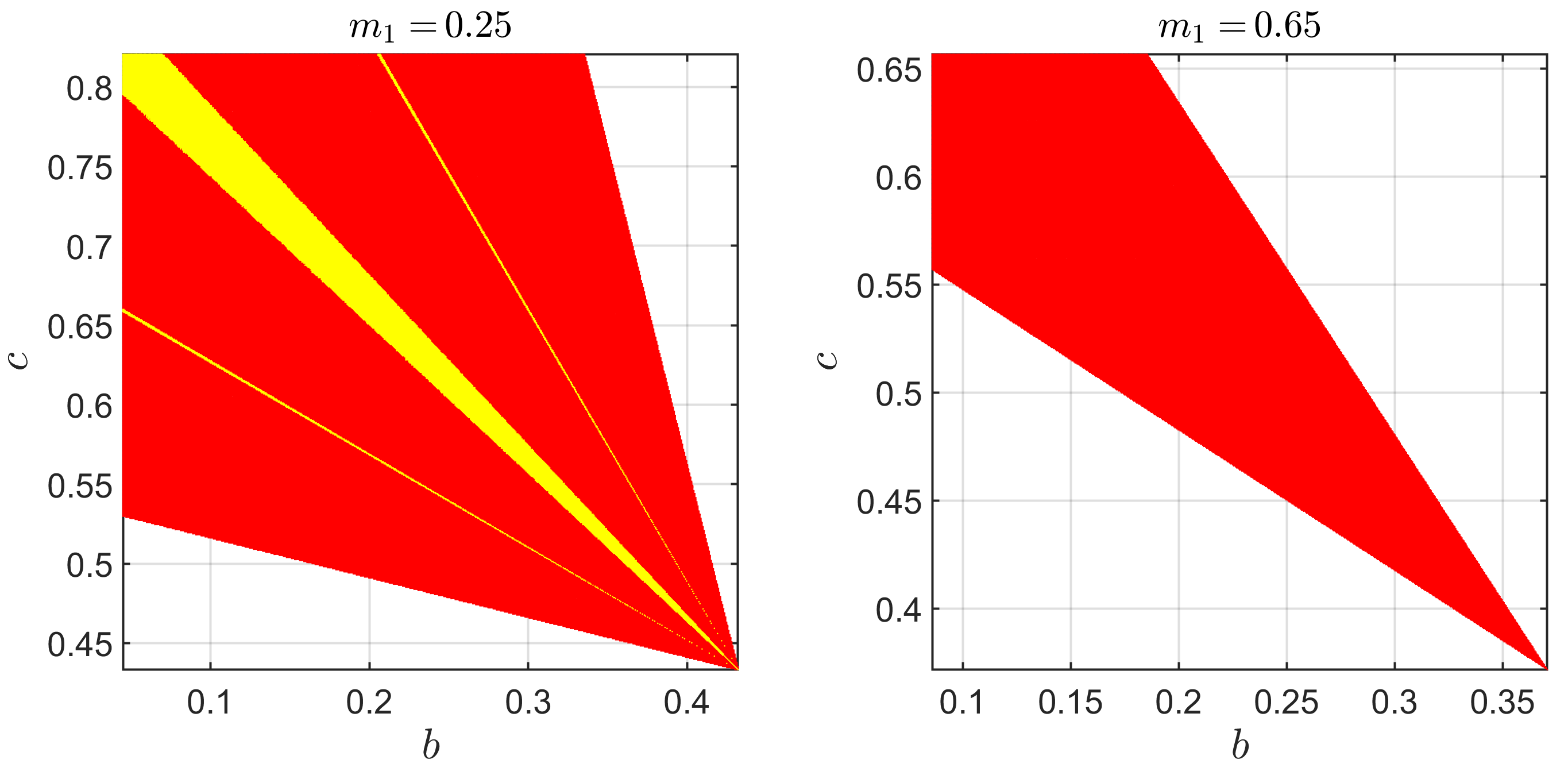}
    \caption{Numerical transitivity test results for 
    the plCNV model with slope parameters 
    ($m_1 = 0.25$, $m_1 = 0.65$). 
    Red points indicate $\beta$-transformations that 
    are numerically transitive. A reduction in 
    numerical transitivity 
    is observed with decreasing slope.}
    \label{fig:transPL}
\end{figure*}
\medskip

In the following figures, each point in the parameter space $(b, c)$ is colour-coded as follows:
\begin{itemize}
    \item \textsc{white:} the map is \emph{not} an expanding Lorenz map (this case is generally not of interest);
    \item \textsc{yellow:} the map is a \emph{numerically nontransitive} expanding Lorenz map;
    \item \textsc{red:} the map is a \emph{numerically transitive} expanding Lorenz map;
    \item \textsc{green:} the map is a \emph{numerically nonLEO} expanding Lorenz map;
    \item \textsc{black:} the map is a \emph{numerically LEO} expanding Lorenz map.
\end{itemize}

For the \emph{plCNV simulations}, 
the model parameters are fixed for both numerical 
transitivity and numerical LEO as
$m_0 = 0.864$, $a = 0.2$, 
and $d = 0.4$, with two choices of
slope values $m_1 = 0.25$ (small) 
and $m_1 = 0.65$ (large), while the variables $(b, c)$ 
are explored over a meshgrid covering the intervals 
defined in Table~\ref{table:plsixcond}. 
For the \emph{nlCNV simulations}, 
parameters are fixed as $a = 0.2$ 
and $d = 0.4$, with two choices of
cubic nonlinearity coefficient,
$\mu = 1$ (low) and $\mu = 2$ (high), and $(b, c)$ 
are similarly explored over a meshgrid satisfying 
the intervals in Table~\ref{table:nlsixcond}.

\subsection{Numerical transitivity simulations}

In the CNV model, the topological transitivity 
and its numerical equivalent refer to 
the ability of the membrane potential dynamics 
to evolve such that, over time, the trajectory starting 
from one value in the state space can approach arbitrarily 
close to any other value within that space. 
In other words, the potential takes on values
that are arbitrarily close to any given value within 
a certain fixed range, which in a sense expresses 
its high flexibility. Such a behavior reflects 
a form of weak chaos, where the system exhibits 
sensitive dependence on initial conditions and 
is capable of producing diverse 
but essentially chaotic firing patterns.

\paragraph{Numerical transitivity for plCNV}
Figure~\ref{fig:transPL} shows the numerical transitivity 
and numerical nontransitivity regions for 
the plCNV model as a function of the invariant interval
endpoints $(b, c)$, using Algorithm~\ref{alg:TransTest}. 
Each point corresponds to a $(b, c)$ pair. 
Red points satisfy both the $\beta$-transformation and
numerical transitivity conditions; yellow points satisfy only 
the $\beta$-transformation condition. 
At first glance, we can see that with a large slope, 
we have numerical transitivity across the entire set 
of acceptable parameters (right panel). 
Whereas with a smaller slope, areas of 
numerical nontransitivity appear (left panel).
It seems that as the slope $m_1$
increases from $0.25$ to $0.65$, the numerically
nontransitive yellow region shrinks to the empty set, 
indicating an increase of dynamical flexibility. 
This suggests that higher slopes
amplify the system's ability to transition between 
different voltage states, increasing the prospect 
for complex and chaotic neuronal behavior.
However, observe that
if a yellow set of nontransitive parameters occurs, 
it occupies consolidated and nondispersed areas. 
Moreover, regular patterns of intertwining transitive 
and nontransitive areas are clearly visible.

\paragraph{Numerical transitivity for nlCNV}

\begin{figure*}[!htb]
    \centering
    \includegraphics[scale=0.21, trim=5mm 5mm 0mm 0mm]{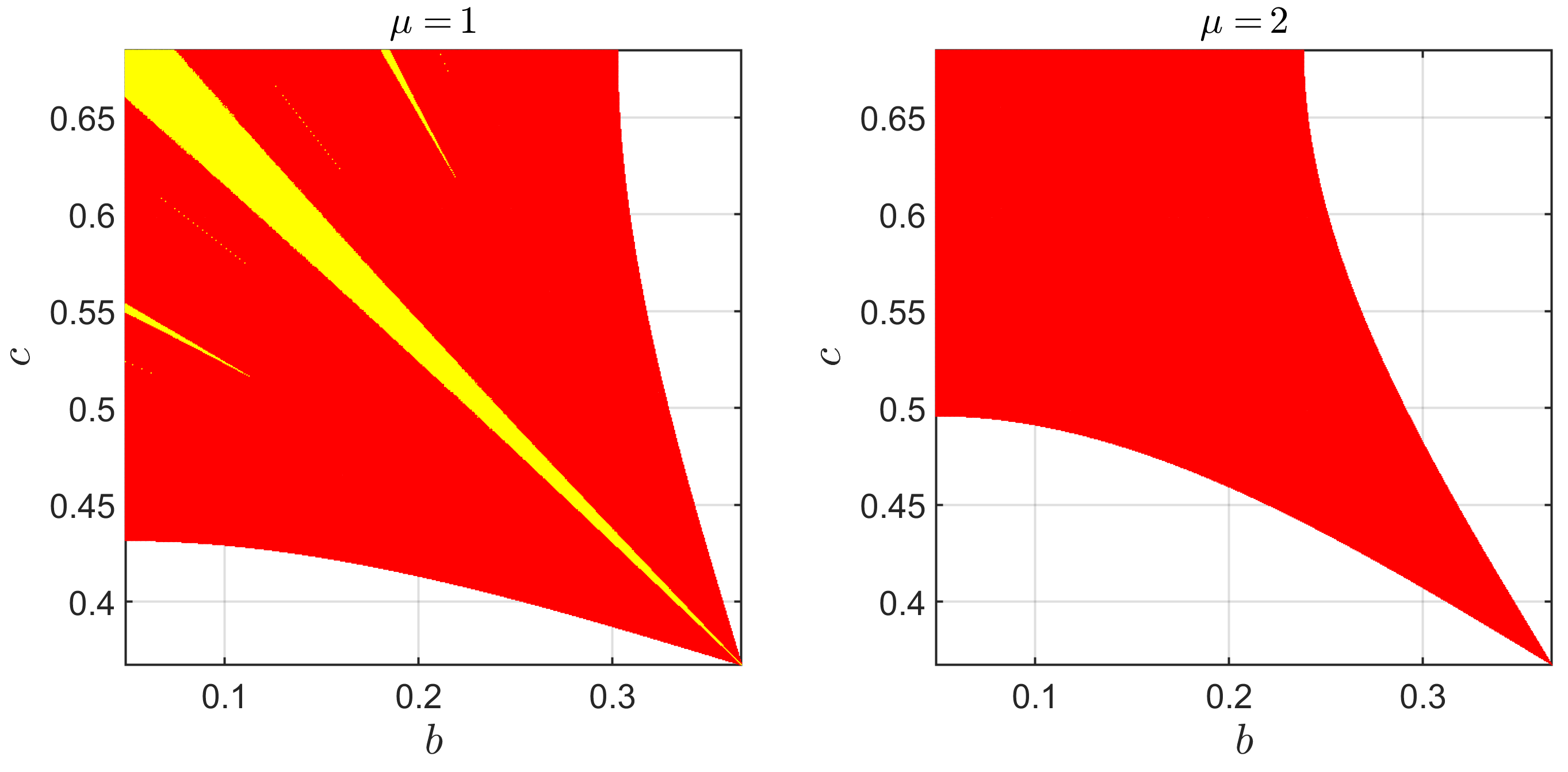}
    \caption{Numerical transitivity test results for the nlCNV model under two different cubic nonlinearity coefficients ($\mu = 1$ and $\mu = 2$). Higher nonlinearity coefficient leads to a complete transitive region, corresponding 
    to more complex and irregular dynamics.}
    \label{fig:transNL}
\end{figure*}

First of all, let us note that in Figure~\ref{fig:transNL}
(left panel) the patterns of 
overlap between transitive and nontransitive areas 
in the set of parameters characteristic of the piecewise 
linear CNV model occur in a very similar way in 
the nonlinear CNV model.
This figure presents numerical transitivity 
simulations for the 1D nlCNV model, where disappearance of
numerical nontransitivity is observed with increasing cubic 
nonlinearity coefficient. 
For $\mu = 2$, the red region takes up the entire 
colored area, indicating richer dynamics and broader 
behavior compared to $\mu = 1$. 
The nlCNV model 
thus shows greater sensitivity and adaptability with 
stronger nonlinearity,
as in the case of larger slope in the plCNV model.

\subsection{Numerical LEO simulations}

\begin{figure*}[!htb]
    \centering
    \includegraphics[scale=0.21, trim=5mm 5mm 0mm 0mm]{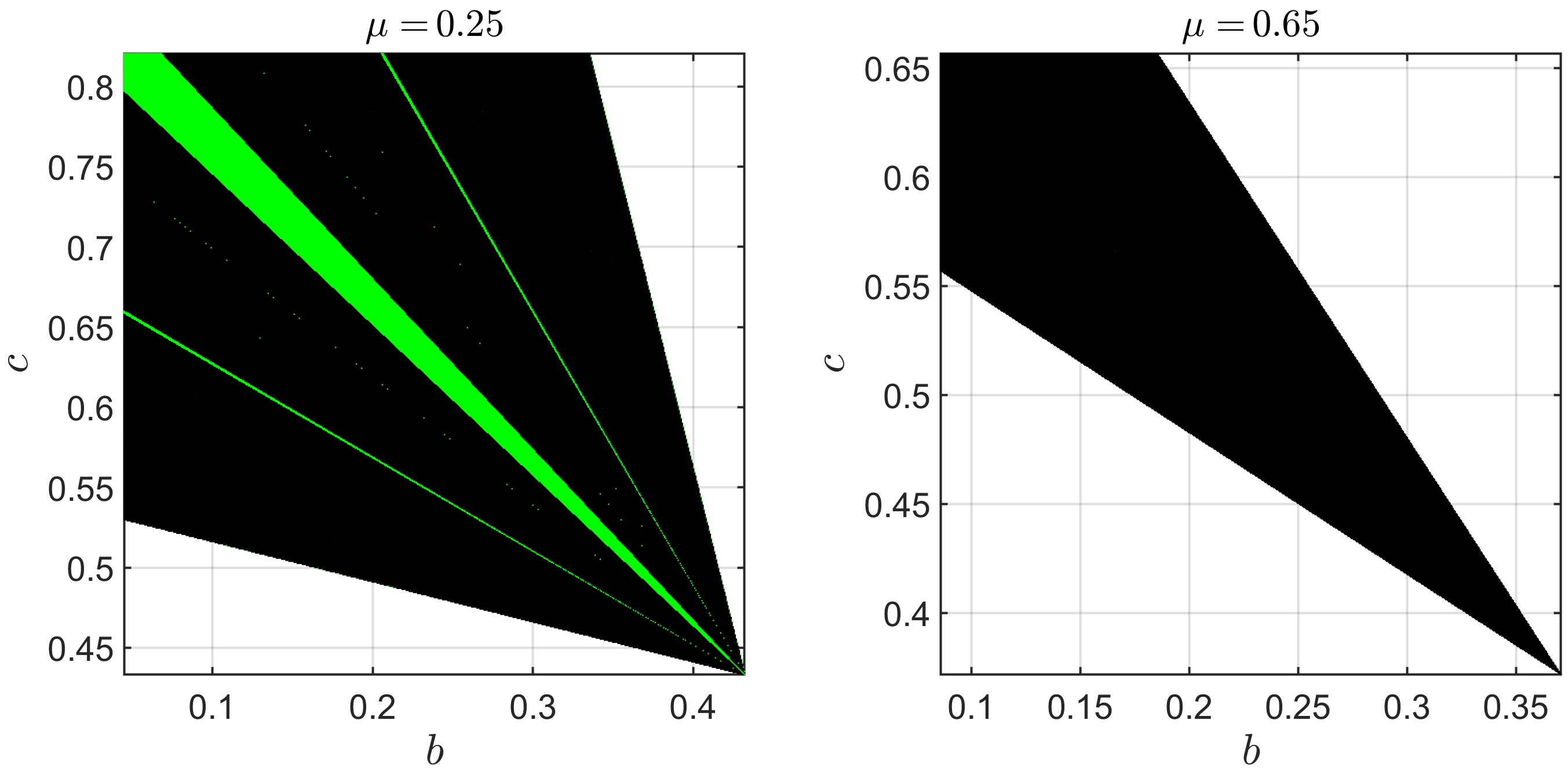}
    \caption{Numerical LEO test for the plCNV model with two slope values ($m_1 = 0.25$ and $m_1 = 0.65$). Black points indicate the numerical LEO and green ones the numerical
    nonLEO. No significant differences can be seen in comparison with Figure~\ref{fig:transPL}.}
    \label{fig:LEO-plCNV}
\end{figure*}

\begin{figure*}[!htb]
    \centering
    \includegraphics[scale=0.21, trim=5mm 5mm 0mm 0mm]{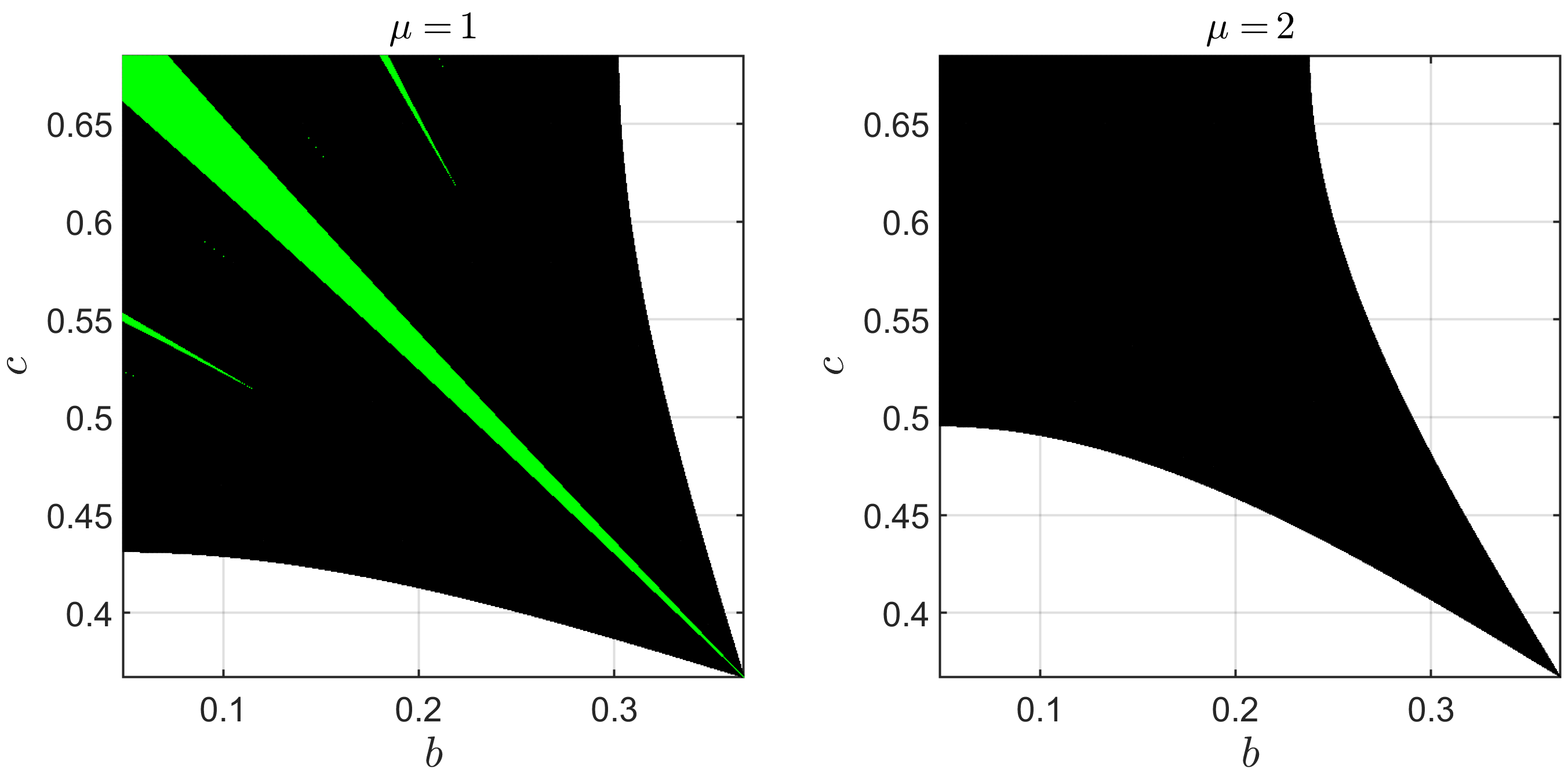}
    \caption{Numerical LEO test results for the nlCNV model with cubic nonlinearity coefficient values $\mu = 1$ and $\mu = 2$. Compare with Figure~\ref{fig:transNL}
    concerning numerical transitivity.}
    \label{fig:LEO-nlCNV}
\end{figure*}

In this subsection we present the results
of numerical LEO simulations for two
versions of the 1D CNV model. Recall that
formally the topological LEO property is stronger 
than the topological transitivity. The topological LEO 
ensures that, starting from any small region, 
the system’s dynamics will eventually 
cover the entire invariant interval. 
Roughly speaking, in neurons, 
the topological LEO behavior represents maximum
unpredictability in firing patterns, 
making it a marker of strong chaos 
and complex neural coding.
We expected that the numerical LEO
would also be a stronger (more restrictive) property
than the numerical transitivity. However, 
a quick comparison of the results of the numerical
transitivity (Figure~\ref{fig:transPL} 
and~\ref{fig:transNL})
and numerical LEO 
(Figure~\ref{fig:LEO-plCNV} and~\ref{fig:LEO-nlCNV})
simulations suggests quite the opposite. 
Namely, it suggests that there are no significant 
differences between the results of both tests.
It seems that the differences 
only concern a few individual points.

\paragraph{Numerical LEO for plCNV}

Figure~\ref{fig:LEO-plCNV} shows the results of
the numerical LEO test for plCNV. 
Figure~\ref{fig:LEO-plCNV} seems to be
almost identical to Figure~\ref{fig:transPL}.
This means that almost always in the space
of parameters the test of numerical transitivity
and the test of numerical LEO give the same result, i.e.,
except for a small number of points the numerical
transitivity test is positive (negative)
if and only if the numerical
LEO test is positive (negative).
In consequence, since the numerical transitivity test
is simpler and much more computationally efficient,
it seems that, in general, there is no point 
in conducting a separate numerical LEO test, 
and it can be completely replaced 
by a numerical transitivity test.

\paragraph{Numerical LEO for nlCNV}

Figure~\ref{fig:LEO-nlCNV} confirms that the same occurs 
in the nonlinear case. Namely, again this figure 
(numerical LEO) seems to be almost identical to 
Figure~\ref{fig:transNL} (numerical transitivity).
However, in a few isolated cases (probably borderline cases), 
the tests give different results.

\medskip
\noindent
To summarize this section, the results from both the plCNV and nlCNV 
models demonstrate how changes in structural or nonlinear parameters 
directly influence the system's capacity for numerical 
transitive and LEO dynamics. These transitions show how sensitive 
the CNV model is to small changes. Such flexibility is essential 
for modeling the diverse firing patterns observed in real neurons.
Both numerical transitivity and LEO reflect chaos and support 
the idea that high dynamical variability can enhance information 
encoding and robustness in neural signal 
processing~\cite{rabinovich2006dynamical}.
\section{Conclusion and discussion}\label{sec:discussion}

Let us briefly summarize our results and discuss some
further perspectives. Firstly, our previous and current 
studies demonstrate that for expanding Lorenz maps, 
topological transitivity is equivalent to chaos, 
making it a primary characteristic of complex dynamics. 
Moreover, the test for numerical transitivity 
was found to be both accurate and efficient in approximating 
topological transitivity, while numerical LEO provided 
a corresponding relationship with topological LEO. However, 
there is one significant difference between these two
numerical tests: computational performance. 
We found that numerical transitivity greatly outperforms 
numerical LEO, yet in nearly all cases, both methods yield 
almost the same results. This close correspondence implies 
that, in many situations, we can substitute the less complicated 
and less computer-intensive test for numerical transitivity 
in place of the numerical LEO test without sacrificing effectiveness. 

\begin{figure*}[!htb]
    \centering
    \includegraphics[scale=0.21, trim=0mm 0mm 0mm 0mm]{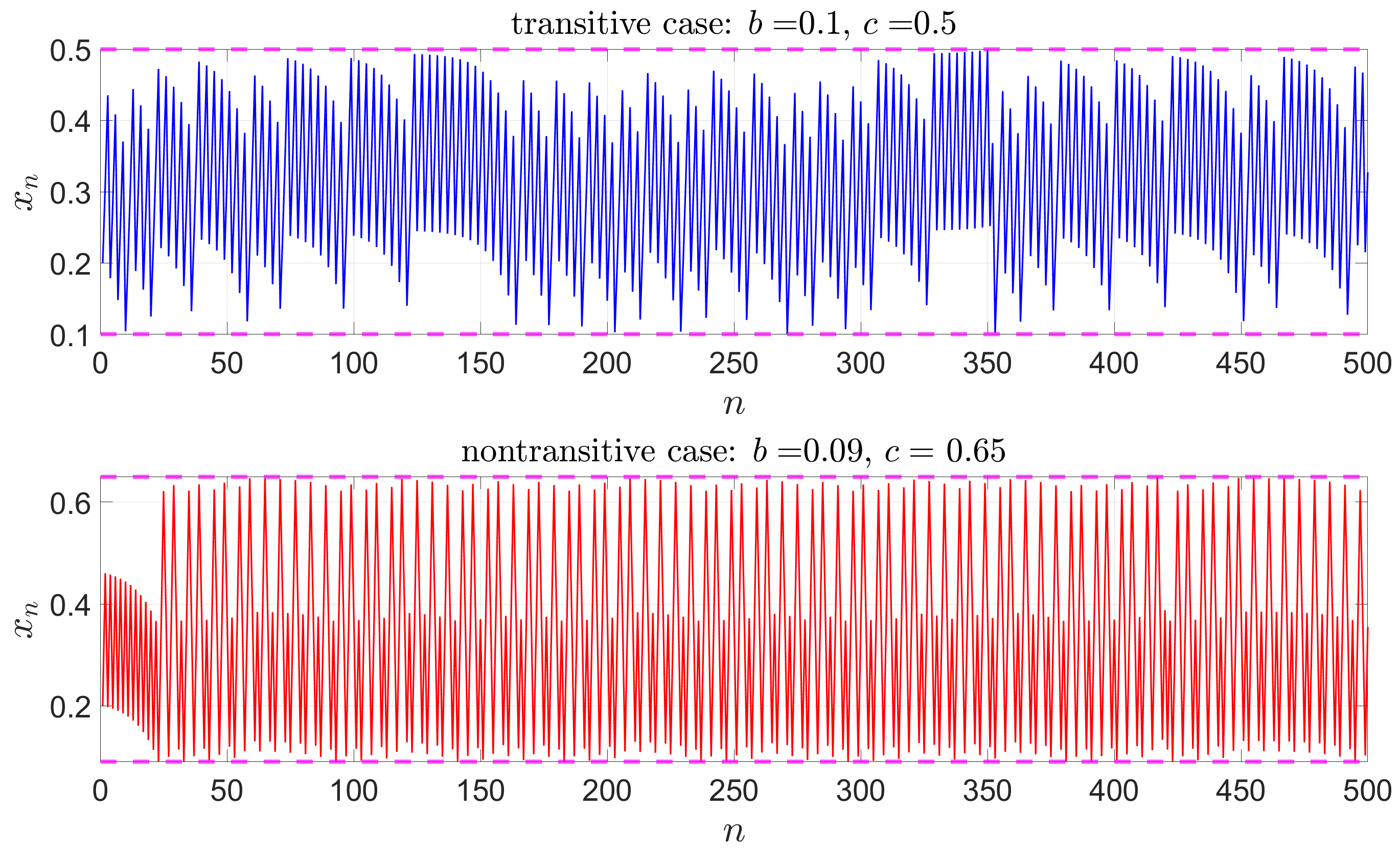}
    \caption{Two different patterns of time series
    of the voltage $x_n$ in the nlCNV model 
    reproducing from the top: 
    bursting (transitive case) and oscillatory
    spiking (nontransitive case). Other 
    parameters: $\mu=0.5$, $a=0.1$, $d=0.37$.}
    \label{fig:timeseries}
\end{figure*}

Next, we explored both numerical tests (transitivity and LEO) 
with respect to the classical family of $\beta$-transformations
(triangle of parameters) and extended their scope to piecewise linear 
and nonlinear CNV neuron models. By combining theoretical examination 
and numerical experiments, we identified parameter ranges where 
transitivity and LEO dynamics emerge and influence neuronal 
behavior. These findings provide a link between theoretical 
descriptions of dynamical systems and biologically relevant models 
and show how deterministic chaos can emerge in firings 
at the level of neurons. In doing so, our research provides deeper 
insight into how mathematical properties like transitivity 
and LEO are related to their biological counterparts 
in real neuronal dynamics.

\begin{figure*}[!htb]
    \centering
    \includegraphics[scale=0.21, trim=5mm 5mm 0mm 0mm]{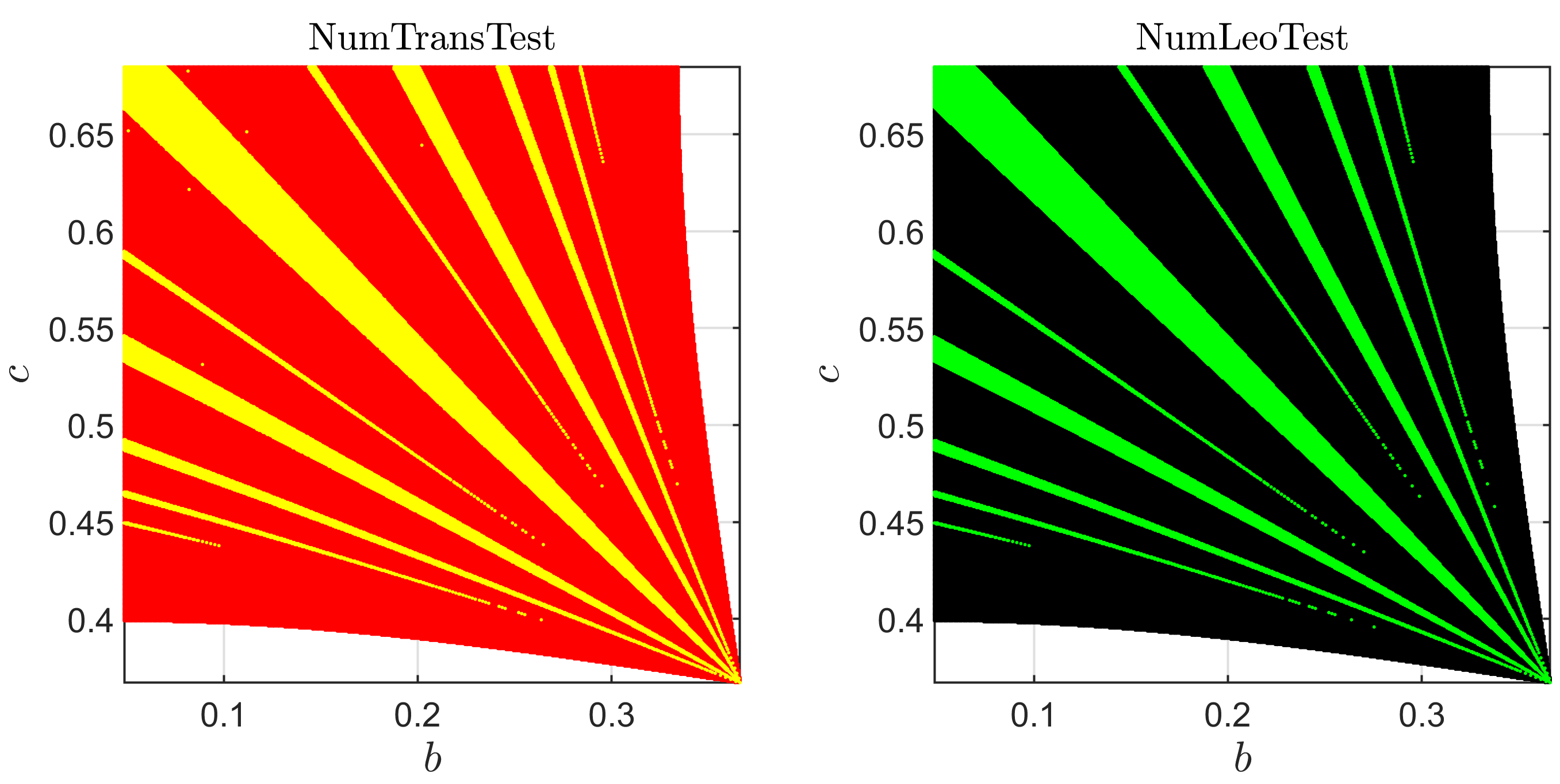}
    \caption{Comparison of results of numerical transitivity (left)
    and numerical LEO (right) tests in the nlCNV
    model. Parameters: $\mu=0.5$, $a=0.1$, $d=0.37$. Meshgrid = 600.}
    \label{fig:transvsleomu}
\end{figure*}

Note that the 1D nlCNV neuron model can
display a variety of different behaviors (firing patterns),
which is important for reproducing the
general spike train patterns such as tonic spiking, 
chaotic bursting, subthreshold oscillations
and others. 
Some of these firing patterns are
depicted in Figure~\ref{fig:timeseries},
where we see \emph{chaotic bursting} 
appearing in the \emph{transitive} regime (top) 
and \emph{oscillatory spiking} 
appearing in the \emph{nontransitive} regime (bottom).
Observe that the top pattern (corresponding to numerically 
transitive parameters) is much more irregular (chaotic)
and tightly fills the range of possible voltages 
with its values, while the bottom pattern 
(corresponding to numerically 
nontransitive parameters) is more regular (almost periodic) 
and bypasses large portions of the range 
of possible voltages. This suggests the possibility of 
quick preliminary classification of firing patterns,
which occur in 1D neuron models,
based on the results of the numerical transitivity test.

\begin{figure}[!htb]
    \centering
    \includegraphics[scale=0.185, trim=0mm 0mm 0mm 0mm]{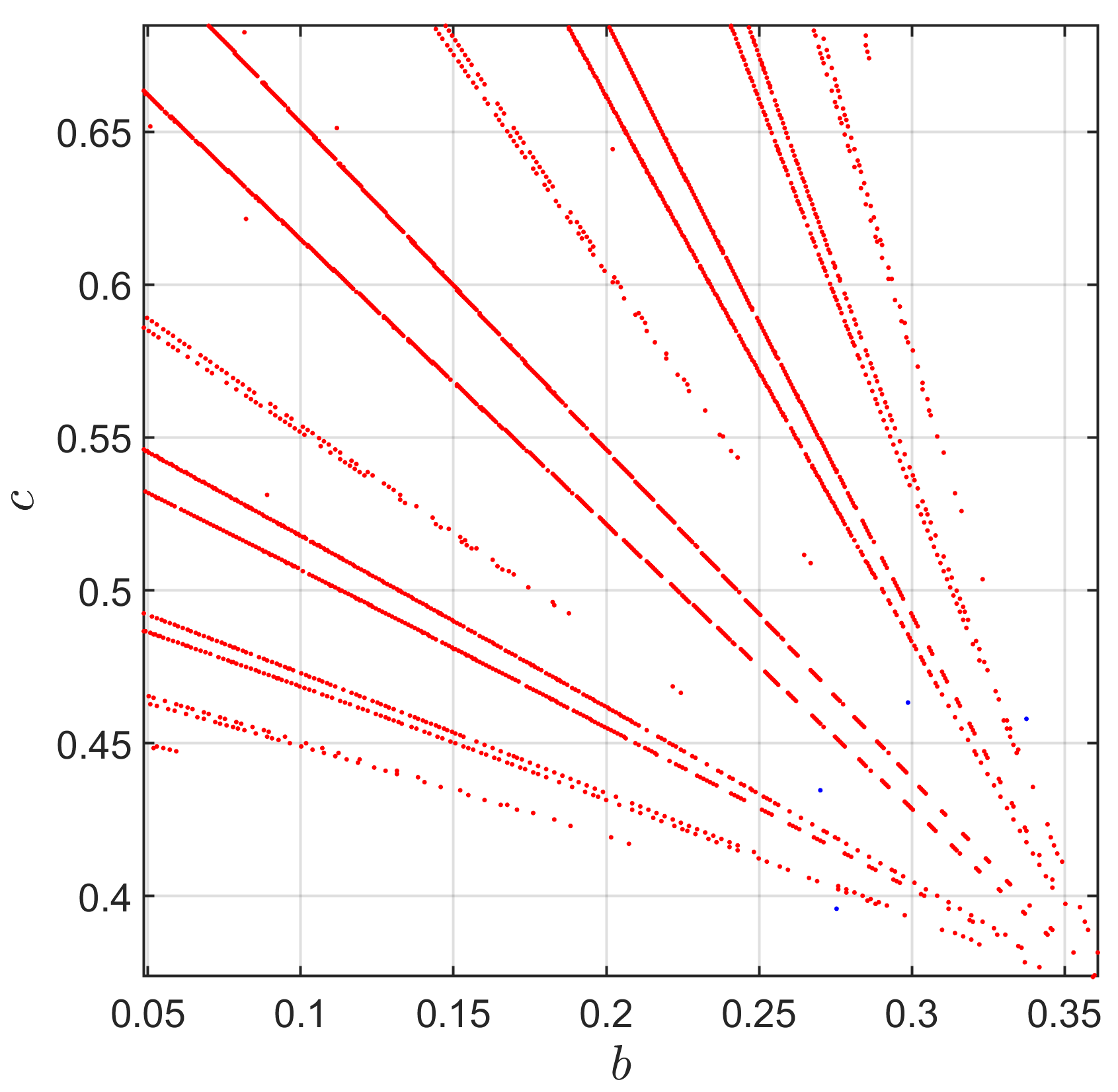}
    \caption{Differences between numerical transitivity 
    and numerical LEO test results  
    in the nlCNV model
    in the $b$–$c$ parameter plane.
    Blue points = transitivity and not LEO
    (there are only four such points),
    red points = LEO and not transitivity. Parameters as
    in Figure~\ref{fig:transvsleomu}.
    Meshgrid = 600.}
    \label{fig:transvsleo}
\end{figure}

In conclusion, both numerical transitivity and numerical LEO 
are effective and reliable tools for analyzing chaotic dynamics 
in expanding Lorenz maps. Their robustness makes them particularly 
valuable for exploring real-world systems modeled within 
this framework, such as the CNV neuron model. 
Among the two, numerical transitivity is a computationally optimal 
approach, offering a straightforward test for chaos that is applicable
across a wide range of parameter settings. Since the results from both
methods are nearly identical in most scenarios
(as shown in Figures~\ref{fig:transvsleomu}
and~\ref{fig:transvsleo}), 
numerical transitivity can generally replace 
the more resource-intensive numerical LEO test. 
These findings strengthen the connection between 
the theoretical properties of dynamical systems and 
their biological applications, highlighting the usefulness 
of simple but powerful numerical approaches in studying 
chaotic behavior in neuronal models.

In comparing our results with other approaches, it is relevant 
to mention the 0–1 test for chaos introduced by Gottwald 
and Melbourne~\cite{GottwaldMelbourne2004, GottwaldMelbourne2009}. 
Due to its simplicity and versatility in applicability to a range 
of dynamical systems, its universal adoption is traceable to 
the fact that it produces a binary outcome indicating either 
presence or absence of chaotic dynamics. However, unlike the 
numerical transitivity and LEO tests, whose outcomes 
are necessarily coupled with the topological properties 
of interval maps, the 0–1 test does not capture finer structural 
differences such as transitivity or LEO. 
This further emphasizes the advantage of our approach in 
correlating numerical outcomes with sound theoretical constructs 
in dynamical systems. 

At the same time, it is to be noted that 
efficiently checking for topological mixing is still yet 
to be achieved by any comparable effort, 
probably due to its complex definition.  
Though topological mixing is a condition stronger than 
topological transitivity and weaker than topological LEO, 
it still lacks a reliable numerical algorithm for detecting it 
see~\cite{spiking_neurons, lorenz_attractor_as_mixing}. 
The establishment of such a computational algorithm 
seems to be quite challenging,
but it would increase our ability 
to classify and understand chaotic motion 
in theoretical and application contexts alike.
\FloatBarrier

\section*{Declaration of competing interest}
The authors declare that they have no known competing financial interests 
or personal relationships that could have appeared to influence 
the work reported in this paper.


\section*{Data availability}
No data was used for the research described in the article.
\FloatBarrier
\bibliographystyle{elsarticle-num-names}
\bibliography{Arxiv}
\end{document}